\documentclass[11pt,a4paper]{article}

\usepackage{amsmath}     
\usepackage{dsfont}		
\usepackage{amsfonts}		
\usepackage{amsthm}
\usepackage{graphicx}
\usepackage{amssymb}

\usepackage{latexsym}
\usepackage{euscript,makeidx,color,mathrsfs}
\usepackage{enumerate}
\usepackage{tabu}
\usepackage[thicklines]{cancel} 

\usepackage{tikz}  
\usetikzlibrary{arrows.meta}
\usetikzlibrary{positioning} 
\usetikzlibrary{graphs}
\usepackage{xcolor} 
\usetikzlibrary {decorations.pathmorphing} 

\usepackage[colorlinks,linkcolor=blue,anchorcolor=green,citecolor=red]{hyperref}

\def\sqr#1#2{{\vcenter{\vbox{\hrule height.#2pt
              \hbox{\vrule width.#2pt height#1pt \kern#1pt \vrule width.#2pt}
              \hrule height.#2pt}}}}
\def\5n{\negthinspace \negthinspace \negthinspace \negthinspace \negthinspace }
\def\4n{\negthinspace \negthinspace \negthinspace \negthinspace }
\def\3n{\negthinspace \negthinspace \negthinspace }
\def\2n{\negthinspace \negthinspace }
\def\1n{\negthinspace }

\def\dbF{\mathbb{F}}

\def\dbH{\mathbb{H}}

\def\dbN{\mathbb{N}}

\def\dbP{\mathbb{P}}

\def\dbR{\mathbb{R}}
\def\dbS{\mathbb{S}}

\def\={\buildrel \triangle \over =}

\def\ds{\displaystyle}

\def\b{\beta}

\def\d{\delta}
\def\e{\varepsilon}

\def\l{\lambda}

\def\f{\varphi}

\def\Si{\Sigma}

\def\O{\Omega}
\def\mf{\mathcal{F}}
\def\me{\mathbb{E}}

\def\rd{\,\mathrm d}  
\def\bal{\begin{aligned}}
\def\eal{\end{aligned}}

\def\cC{{\cal C}}

\def\cF{{\cal F}}

\def\cJ{{\cal J}}

\def\cT{{\cal T}}
\def\cU{{\cal U}}

\def\no{\noindent}

\def\ms{\medskip}

\def\q{\quad}
\def\qq{\qquad}

\def\liminf{\mathop{\underline{\rm lim}}}

\def\lt{\left}
\def\rt{\right}
\def\lan{\langle}
\def\ran{\rangle}

\def\rf{\eqref}

\def\h{\widehat}
\def\wt{\widetilde}

\def\cd{\cdot}
\def\cds{\cdots}

\def\ae{\hbox{\rm a.e.}}

\def\deq{\triangleq}

\def\({\Big (}
\def\){\Big )}
\def\[{\Big[}
\def\]{\Big]}

\def\bde{\begin{definition}\label}
\def\ede{\end{definition}}
\def\be{\begin{equation}}
\def\bel{\begin{equation}\label}
\def\ee{\end{equation}}
\def\beq{\begin{equation*}\begin{aligned}}
\def\eeq{\end{aligned}\end{equation*}}
\def\bt{\begin{theorem}\label}
\def\et{\end{theorem}}
\def\bc{\begin{corollary}\label}
\def\ec{\end{corollary}}
\def\bl{\begin{lemma}\label}
\def\el{\end{lemma}}
\def\bp{\begin{proposition}\label}
\def\ep{\end{proposition}}
\def\bas{\begin{assumption}\label}
\def\eas{\end{assumption}}
\def\br{\begin{remark}\label}
\def\er{\end{remark}}
\def\bex{\begin{example}\label}
\def\ex{\end{example}}
\def\ba{\begin{array}}
\def\ea{\end{array}}
\def\ed{\end{document}}

\def\square#1{\vbox{\hrule\hbox{\vrule height#1%
     \kern#1\vrule}\hrule}}
\def\rectangle#1#2{\vbox{\hrule\hbox{\vrule height#1%
     \kern#2\vrule}\hrule}}

\font\tenbb=msbm10 \font\sevenbb=msbm7 \font\fivebb=msbm5

\newfam\bbfam
\scriptscriptfont\bbfam=\fivebb \textfont\bbfam=\tenbb
\scriptfont\bbfam=\sevenbb

\newtheorem{theorem}{\hskip 1.3em Theorem}[section]
\newtheorem{definition}[theorem]{\hskip 1.3em Definition}
\newtheorem{proposition}[theorem]{\hskip 1.3em Proposition}
\newtheorem{corollary}[theorem]{\hskip 1.3em Corollary}
\newtheorem{lemma}[theorem]{\hskip 1.3em Lemma}
\newtheorem{remark}[theorem]{\hskip 1.3em Remark}
\newtheorem{example}[theorem]{\hskip 1.3em Example}

\newtheorem{assumption}[theorem]{\hskip 1.3em Assumption}

\allowdisplaybreaks

\makeatletter
   
   \@addtoreset{equation}{section}
\makeatother

\begin{document}

\title{Equivalence between time and norm optimal control problems of stochastic differential equations
}

\author{
Yanqing Wang\thanks{
School of Mathematics and Statistics, Southwest University, Chongqing 400715, China.  {\small\it
e-mail:} {\small\tt yqwang@swu.edu.cn}. \ms}}

\date{Nov. 29, 2024}
\maketitle

\begin{abstract}
In this paper, we establish the equivalence between minimal time and minimal norm control problems for controllable stochastic differential equations (SDEs). The core of this equivalence lies in proving the continuity of the minimal norm function. Unlike the deterministic setting, where the time translation technique is effective, this approach is not applicable due to the adapted nature of the admissible controls in SDEs. To address this challenge, we construct the minimal norm control explicitly, thereby achieving the desired continuity. This constructive approach not only proves theoretical equivalence but also facilitates the design of efficient numerical schemes for minimal norm and time controls.

\end{abstract}

\ms

\no\bf Keywords: \rm 
equivalence, minimal norm control, minimal time control, stochastic differential equation

\ms

\no\bf AMS 2020 subject classification: \rm  
49N10, 93B52
 



\section{Introduction}\label{sec-intro}

Let $(\Omega, {\mathcal F}, {\mathbb P})$ be a complete probability space on which a 
standard one-dimensional Brownian motion $W=\{W(t)\,|\,t\geq 0\}$ is defined, 
with $\dbF=\{\cF_t\}_{t\geq 0}$ being its natural filtration augmented by all the $\dbP$-null sets in $\cF$.
For fixed $T\in (0,\infty)$ and a given Hilbert space $(\dbH, \lan\cd,\cd\ran_{\dbH})$, we denote by $L^2_{\mf_t}(\O;\dbH)$ for $t\in [0,T]$ the Hilbert space consisting of all $\mf_t$-measurable, $\dbH$-valued random variables $\xi$
such that $\me\|\xi\|_\dbH^2<\infty$; by $L^2_\dbF(0,T;\dbH)$ the Hilbert space consisting of all $\dbH$-valued, $\dbF$-adapted stochastic processes $\f(\cd)$ such that $\me\int_0^T\|\f(t)\|_\dbH^2\rd t<\infty$; by $L^2_\dbF(\O;C([0,T];\dbH))$ the Banach space consisting of continuous processes $\f(\cd)\in L^2_\dbF(0,T;\dbH)$
such that $\me\big[\sup_{t\in[0,T]}\|\f(t)\|_\dbH^2\big]<\infty$. In a similar way, one can define $L^2_\dbF(\dbR^+;\dbH)$.

This paper is devoted to study the equivalence of two optimal control problems --- minimal {\em time} and {\em norm} control problems --- subject to the following linear controlled stochastic differential equation (SDE, for short) 
\bel{sys-1}
\bal
&{\mathrm d} x(t)=\big[Ax(t)+Bu(t)\big]\rd t+\big[Cx(t)+Du(t)\big]\rd W(t) \qq t\geq 0\,,\\
\eal
\ee
where $A\,,C\in \dbR^{n\times n}$, $B\,,D\in \dbR^{n\times m}$ for $n,m\in \dbN^+$, and $u(\cd)$ is a
 control variable taken from 
$L^2_\dbF(\dbR^+;\dbR^m)$.
For a given initial condition $x(0)=x_0\in\dbR^n$, we will denote by $x\big(\cd; x_0,u(\cd)\big)$ 
the solution to SDEs \rf{sys-1}.

Intuitively, the minimal time control problem focuses on reaching a desired target in the shortest possible time, while the minimal norm control problem aims to achieve the target with the least amount of energy or control effort.
The solvability of these two problems for systems governed by SDEs, and 
the exact controllability of linear controlled SDEs, along with
the validity of 
 observability inequalities for backward SDEs (BSDEs, for short) are four notable problems in stochastic control theory. In \cite{Wang-Yang-Yong-Yu17}, the equivalence of the solvability of minimal norm problems, controllability of SDEs and the validity of observability inequalities for BSDEs was established. 
Additionally, there exist references addressing one or more of these topics, such as \cite{Peng94, Liu-Peng10,Gashi15, Wang-Zhang15,Dou-Lv19, Wang-Yu20, Lv-Zhang21}.
In the deterministic setting, some works have explored the equivalence between minimal time and norm control problems; see \cite{Wang-Zuazua12, Wang-Xu-Zhang15, Qin-Wang18} and the references therein.
However,
to the best of the author’s knowledge, there are no published works
 on the equivalence of the two problems in the context of SDEs.
Therefore, this work aims to fill this gap.

\ms

If $C=0$, $D=0$ and $u(\cd)\in L^2(\dbR^+;\dbR^m)$, \rf{sys-1} is
reduced to a linear controlled ordinary differential equation (ODE, for short), which will be denoted by $[A, B]$.
The controllability of $[A, B]$ is well-studied. However,
differing from the equivalence of exact, approximate and null controllability for the ODE system $[A, B]$, the controllability of SDEs is more complicated and not yet completely understood. For studies on exact controllability,
see e.g., \cite{Peng94,Liu-Peng10,Lv-Yong-Zhang12,Wang-Yang-Yong-Yu17}; for approximate and null controllability,
see e.g., \cite{Buckdahn-Quincampoix-Tessitore06, Goreac16}. 
Throughout this work, we adopt the following assumption:\\
{\bf (A)} System \rf{sys-1} is exactly controllable, {\em i.e.}, for any $T>0$, $x_0\in\dbR^n$, $\xi\in L^2_{\mf_T}(\Omega;\dbR^n)$, there exists a control
$u(\cd)\in L^2_\dbF(\dbR^+;\dbR^m)$ such that $x\big(T; x_0,u(\cd)\big)=\xi$.

\ms
In the below, we specify minimal {\em time} and {\em norm} control problems governed by SDE \rf{sys-1}.
Suppose that $x_0\in \dbR^n$ is any given initial state and $\xi\in L^2_{\mf_{T_0}}(\Omega;\dbR^n)$ is any fixed target with
$T_0\in \dbR^+$.
It is obvious that if $T_0=0$, $\xi\in \dbR^n$ is a deterministic vector.
Now the minimal time control problem
subject to \rf{sys-1} can be stated as Problem $(TP)^N_{x_0,\xi}$: 
for any $N\geq 0$,
\bel{minimal-time}
T(N,x_0,\xi)\deq \inf_{u(\cd)\in\cU_N}\big\{t\in \dbR^+\,\big|\, x\big(t;x_0,u(\cd)\big)=\xi\big\}\,,
\ee
where the feasible control set $\cU_N$ is defined by
\bel{w1227e1}
\cU_N\deq \big\{ u(\cd)\in L^2_\dbF(\dbR^+,\dbR^m)\,\big|\, \|u(\cd)\|_{L^2_\dbF(\dbR^+;\dbR^m)}\leq N \big\}\,.
\ee
If $\big\{t\in \dbR^+\,|\, x\big(t;x_0,u(\cd)\big)=\xi \,, \mbox{for any } u(\cd)\in \cU_N\big\}=\emptyset$, set
$T(N,x_0,\xi)=\inf \emptyset=\infty$. 

\bde{def-tp}
$T(N,x_0,\xi)$ given in \rf{minimal-time} is called the {\em minimal time} of Problem $(TP)_{x_0,\xi}^N$. If $u(\cd)\in \cU_N$ and there exists some 
$t\in \dbR^+$
such that $x\big(t; x_0, u(\cd)\big)=\xi$, then $u(\cd)$ is called an {\em admissible control}. An admissible control $u(\cd)$ is called 
a {\em minimal time control} if $T(N,x_0,\xi)<\infty$ and $x\big(T(N,x_0,\xi), x_0, u(\cd)\big)=\xi$.
If the restrictions of all minimal time controls over $\big(0, T(N,x_0,\xi)\big)$ are the same, then the minimal time control is 
said to be {\em unique}. 
\ede

By the adaptedness property of solutions to SDEs, we know that for any $t<T_0\,, N\geq 0$ and $u(\cd)\in \cU_N$, 
it holds that $x\big(t; x_0, u(\cd)\big)\in L^2_{\mf_t}(\Omega;\dbR^n)\subsetneq L^2_{\mf_{T_0}}(\Omega;\dbR^n)$. 
Thus, if $L^2_{\mf_{T_0}}(\Omega;\dbR^n)$ is the ``smallest" space where $\xi$ belongs to, {\em i.e.}, 
\bel{w1226e1}
\xi\notin L^2_{\mf_{t}}(\Omega;\dbR^n) \qq \forall\, t\in [0,T_0)\,,
\ee
then
it is impossible to find a control $u(\cd)\in \cU_N$ such that $x(t; x_0, u(\cd))=\xi$ for any $t<T_0$. 
Equivalently, this implies
\beq
T(N,x_0,\xi)\geq T_0 \qq \forall\, N\geq 0\,.
\eeq
For example, if the target $\xi=\e W(T_0)$ with $\e>0$, no matter how small $\e$ is or how great the control capability is, $\xi$ cannot be reached before time $T_0$.

Similar to Problem $(TP)_{x_0,\xi}^N$, we also can define 
the minimal  norm control problem
stated as Problem $(NP)_{x_0,\xi}^T$: for any $T>0$,
\bel{minimal-norm}
N(T,x_0,\xi)\deq \inf\big\{ \|u(\cd)\|_{L^2_\dbF(0,T;\dbR^m)}\,\big|\, x\big(T;x_0,u(\cd)\big)=\xi \big\}\,.
\ee

\bde{def-np}
In Problem $(NP)_{x_0,\xi}^T$, $N(T,x_0,\xi)$ defined in \rf{minimal-norm} is called the {\em minimal norm}. 
If $x\big(T;x_0,u(\cd)\big)=\xi$, then the $L^2_\dbF(0,T;\dbR^m)$-valued $u(\cd)$ is called an {\em admissible control}.
An admissible control $u(\cd)$ is called a
{\em minimal norm control} if $\|u(\cd)\|_{L^2_\dbF(0,T;\dbR^m)}=N(T,x_0,\xi)$.
\ede

\bde{equiv}
Problems $(TP)_{x_0,\xi}^N$ and $(NP)_{x_0,\xi}^T$ are said to be mutually equivalent if the following conditions hold:
\begin{enumerate} 
\item [ \rm (i)] Problems $(TP)_{x_0,\xi}^N$ and $(NP)_{x_0,\xi}^T$ have minimal time and norm controls, respectively;

\item [ \rm (ii)] When $u_*(\cd)$ is a minimal time control of Problem $(TP)_{x_0,\xi}^N$, then $u_*|_{(0,T)}(\cd)$, the restriction of $u_*(\cd)$ to the interval $(0,T)$, is a minimal norm control of Problem $(NP)_{x_0,\xi}^T$;

\item [ \rm (iii)] When $u^*(\cd)$ is a minimal norm control of Problem $(NP)_{x_0,\xi}^T$, then $\wt u^*(\cd)$, defined on $\dbR^+$ by extending $u^*(\cd)$ to zero outside of $(0,T)$, is a minimal time control of Problem $(TP)_{x_0,\xi}^N$.

\end{enumerate}
\ede

Compared to  deterministic minimal time and norm control problems, the stochastic counterparts  enjoy their own characteristics:

{\bf (a)} Construction of the minimal norm control. 
In the deterministic minimal norm control problem for a system presented by $[A, B]$, for any initial state $x_0$ and target $\xi$, the minimal norm control $u^*(\cd)\in L^2(0, T;\dbR^m)$ can be presented by
\beq
u^*(t)=- B^\top e^{A^\top (T-t)}\Big[\int_0^T e^{A s}BB^\top e^{A^\top s}\rd s\Big]^{-1}\big[e^{A T}x_0-\xi\big]\qq \forall\, t\in[0,T]\,.
\eeq
In the stochastic setting, however, the above method cannot be directly used to construct admissible controls. 
Driving the presentation of the minimal norm control in the stochastic context is a non-trivial task, even for specific cases; see e.g., \cite{Bi-Sun-Xiong20}. 
In Section \ref{general-target}, we will adopt the method of Lagrange multipliers to derive the representation of the minimal norm control; see Proposition \ref{w1229p4} for details.

{\bf (b)} Adaptedness property of $x(\cd)$ and $u(\cd)$. 
Intuitively,  for any $t\in\dbR^+$, $g(t)$ is $\mf_t$-measurable, where $g(\cd)$ represents either $x(\cd)$ or $u(\cd)$. Hence, after applying a time translation $t\to t+t_0$ for some $t_0>0$,
$g_0(t)\deq g(t+t_0)$ becomes $\mf_{t+t_0}$-measurable, but may no longer be $\mf_t$-measurable. 
To study the equivalence of Problems $(TP)_{x_0,\xi}^N$ and $(NP)_{x_0,\xi}^T$, one needs continuity of the minimal norm function $N(\cd,x_0, \xi)$. In the deterministic
setting,  time translation is a common technique; see e.g., \cite[Theorem 2.6]{Qin-Wang18}. But in the stochastic setting, this technique fails due to the loss of adaptedness of $g_0(\cd)$.
In this work, we apply the derived representation of the minimal norm control in Proposition \ref{w1229p4} to derive the continuity of the minimal norm function; see Theorem \ref{w1230t1}.

{\bf (c)} Forward vs. backward equations.
In the deterministic ODE setting, by the time translation $s=T-t$, one can derive a backward ODE. Since this translation is invertible, methods to study optimal control problems for ODEs and backward ODEs are essentially the same. But in the stochastic setting, the solution to a BSDE  is a stochastic pair with two components; 
for example,  $\big(y(\cd), Y(\cd)\big)$ is the solution to BSDE \rf{bsde-1}. The second component of the solution to 
the concerned system (from the backward perspective) may play different roles in the same problem --- that could be a component of the control variable from the forward 
perspective, or be a component of the state variable from the backward perspective.
For example, in Problem $(\wt{NP})_{x_0,\xi}^T$ in Section \ref{general-target}, from the backward perspective $z_0(\cd)$ in \rf{sys-02} is a 
part of the state $\big(x(\cd),z_0(\cd)\big)$, while by viewing \rf{sys-02} as a SDE $z_0(\cd)$ is a part of the control $\big(z_0(\cd),v_0(\cd)\big)$.
For further references, see e.g., \cite{Yong-Zhou99, Lim-Zhou01, Wang-Yu20}.

\ms

To deal with the equivalence of Problems $(TP)_{x_0,\xi}^N$ and $(NP)_{x_0,\xi}^T$, we introduce another optimal control problem. 
For any $T>0$, consider the following BSDE
\bel{bsde-1}
\lt\{
\bal
& \rd y(t)=-\big[A^\top y(t)+C^\top Y(t)\big]\rd t+Y(t)\rd W(t)\qq t\in[0,T]\,,\\
&y(T)=\eta\in L^2_{\mf_T}(\Omega;\dbR^n)\,,
\eal
\rt.
\ee
and define a cost functional
\beq
J^T(\eta;x_0,\xi)=\frac 1 2 \big\| B^\top y(\cd)+D^\top Y(\cd) \big\|_{L^2_\dbF(0,T;\dbR^m)}^2+\lan x_0, y(0;\eta,T) \ran-\me\lan \xi,\eta \ran\,.
\eeq
We now state the optimal control problem associated with BSDE \rf{bsde-1} and $J^T(\cd)$, which we denote as Problem $(OP)^T_{x_0,\xi}$:
\beq
V(T,x_0,\xi)\deq \inf_{\eta\in L^2_{\mf_T}(\Omega;\dbR^n)} J^T(\eta; x_0,\xi)\,.
\eeq

The following result,  from \cite[Theorems 4.1, 4.5, 5.1]{Wang-Yang-Yong-Yu17}, is on the solvability and equivalence of Problems $(NP)_{x_0,\xi}^T$ and $(OP)_{x_0,\xi}^T$.
\bl{w727l1}
For any $T>0$, $x_0\in\dbR^n$ and $\xi\in L^2_{\mf_{T_0}}(\Omega;\dbR^n)$ with $T_0\in [0, T]$, the following statements are equivalent:
\begin{enumerate}
\item [{\rm (i)}] System \rf{sys-1} is exactly controllable  on $[0,T]$;

\item [{\rm (ii)}] The following observability inequality holds
\bel{exact-ob}
\me\big[\|\eta\|^2\big]\leq \cC(T)\,\me\int_0^T \big\|B^\top y(t)+D^\top Y(t)\big\|^2\rd t \qq \forall\,\eta\in L^2_{\mf_T}(\Omega;\dbR^n)\,,
\ee
where $\big(y(\cd),Y(\cd)\big)$ solves BSDE \rf{bsde-1}, and $\cC(T)$ is a positive constant depending on $T$;

\item [{\rm (iii)}] Problem $(NP)_{x_0,\xi}^T$ has a unique minimal norm control $u^*(\cd)\in L^2_\dbF(0,T;\dbR^m)$;

\item [{\rm (iv)}]  Problem $(OP)_{x_0,\xi}^T$ admits a unique solution $\eta^*\in L^2_{\mf_T}(\Omega;\dbR^n)$. 

\end{enumerate}
Moreover, the minimal norm control $u^*(\cd)$ to Problem $(NP)_{x_0,\xi}^T$ has a representation
\bel{w1226e3}
u^*(\cd)=B^\top y^*(\cd)+D^\top Y^*(\cd)\,,
\ee
where $\big(y^*(\cd),Y^*(\cd)\big)$ solves BSDE \rf{bsde-1} with $y^*(T)=\eta^*$, and
\beq
V(T,x_0,\xi)=-\frac 1 2 N^2(T,x_0,\xi)\,.
\eeq

\el

The rest of this paper is organized as follows. In Section \ref{general-target}, we construct
the optimal tuple for an auxiliary optimal control problem $(\wt {NP})_{x_0,\xi}^T$, whose cost functional achieves a minimum equal to the minimal norm of Problem $(NP)_{x_0,\xi}^T$. By establishing the continuity of this minimum with respect to the time variable $T$, we demonstrate the equivalence of Problems $(TP)_{x_0,\xi}^N$ and $(NP)_{x_0,\xi}^T$. In Section \ref{zero-target}, we focus on the minimal norm function $N(T,x_0,\xi)$ for the stable target $\xi=0$ and derive its strict monotonicity, which provides a basis for calculating the minimal time control.

\section{Equivalence of minimal time and norm control problems with a general target}\label{general-target}

In this section, we primarily address the equivalence between minimal time and norm control problems with a general target. In the deterministic setting, previous works such as \cite{Wang-Xu-Zhang15, Qin-Wang18} focused on optimal control problems with the stable target $\xi=0$, and utilized an existence-proof method. A key technique in those studies was the time translation method, used to establish the continuity of the minimal norm function to the time variable $T$. However, this technique is not applicable to our problems involving a general target. To address this challenge, we explicitly construct the minimal norm control, which enables us to establish the desired continuity of the minimal norm function. This continuity is essential for deriving the equivalence between minimal time and norm control problems.

For any given $T>0$ and $x_0\in\dbR^n$, the state to system \rf{sys-1} can be transferred to $x(T; x_0,0)$ at time $T$. Hence if $\xi=x(T;x_0,0)$, it is evident that the unique minimal norm of Problem $(NP)_{x_0,\xi}^T$ is $u^*(\cd)=0$. The following result formalizes this observation.

\bl{w1213l1}
Under the assumption {\bf (A)}, the following statements are equivalent:\\
{\rm (i)} $\Psi(T)x_0\neq \xi$, where $\Psi(\cd)$ solves the following $\dbR^{n\times n}$-valued SDE
\bel{sde-Psi}
\lt\{
\bal
&\mathrm{d} \Psi(t)=A\Psi(t)\rd t+C\Psi(t)\rd W(t) \qq t\in [0,T]\,,\\
&\Psi(0)=I_n\,.
\eal
\rt.
\ee
{\rm (ii)} $u^*(\cd)\neq 0$, where $u^*(\cd)$ is the minimal norm control of Problem $(NP)_{x_0,\xi}^T$.\\
{\rm (iii)} $\eta^*\neq 0$, where $\eta^*$ is the minimizer of Problem $(OP)_{x_0,\xi}^T$.
%
%
%
\el

\begin{proof}

{\bf (i)$\Rightarrow$(ii)}. 
By contradiction, we suppose that the minimal norm control $u^*(\cd)=0$. Then $x^*\big(\cd;x_0,u^*(\cd)\big)$, the corresponding optimal state, satisfies $x^*\big(t;x_0,u^*(\cd)\big)=\Psi(t)x_0$ for all $t\in[0,T]$.
Note that $\Psi(\cd)$ is invertible; see \cite[Theorem 6.14, Chapter 1]{Yong-Zhou99}.
Therefore the condition that $\Psi(T)x_0\neq \xi$ leads to
\beq
x^*\big(T,x_0,u^*(\cd)\big)
=\Psi(T)x_0\neq \xi\,,
\eeq
which contradicts the fact that $u^*(\cd)$ can be used to transfer the state to $\xi$ at $T$.
That settles the assertion (ii).

\ms
{\bf (ii)$\Rightarrow$(iii)}. 
By contradiction, we suppose that $\eta^*=0$. Then the unique solution of BSDE \rf{bsde-1} $\big(y^*(\cd),Y^*(\cd)\big)=0$, 
and by Lemma \ref{w727l1},
we know that the minimal norm control of Problem $(NP)_{x_0,\xi}^T$
\beq
u^*(\cd)=B^\top y^*(\cd)+D^\top Y^*(\cd)=0\,,
\eeq
that contradicts the statement (ii). Hence, (iii) holds. 

\ms

{\bf (iii)$\Rightarrow$(i)}. 
Firstly, observability inequality \rf{exact-ob} and the representation of the minimal norm control $u^*(\cd)$ in \rf{w1226e3} lead to
\bel{w1226e7}
0< \me\big[\|\eta^*\|^2\big]\leq \cC(T)\, \me\int_0^T\big\|B^\top y^*(t)+D^\top Y^*(t)\big\|^2\rd t
= \cC(T)\,\me\int_0^T \|u^*(t)\|^2\rd t\,,
\ee
which means that the minimal norm control $u^*(\cd)\neq 0$. 
By contradiction, suppose that $\Psi(T)x_0=\xi$. Then it is obvious that the minimal norm control of Problem $(NP)_{x_0,\xi}^T$ is
 $u^*(\cd)=0$, which contradicts \rf{w1226e7}.
That completes the proof.
%
%
%
%
%
%
%
%
\end{proof}

%

For any given $T>0$, $x_0\in\dbR^n$ and $\xi\in L^2_{\mf_T}(\O;\dbR^n)$ satisfying \rf{w1226e1}, 
define
\bel{w1225e14}
\cJ_N\deq \big\{t\in [T_0,\infty)\,\big|\, N(t,x_0,\xi)\leq N\big\} \qq\forall\, N\geq 0\,.
\ee
The following property shows the relationship between $\cJ_\cd$ and the minimal time control $T(\cd,x_0,\xi)$. 

\bl{w1225l3} 
For any $N\geq 0$, it holds that
\beq
T(N,x_0,\xi)=\inf \cJ_N \,.
\eeq
\el

\begin{proof}
The proof draws on ideas from \cite[Theorem 2.4]{Qin-Wang18}, and we present it by considering two separate cases.

\no{\bf Case I.} $\cJ_N=\emptyset$

In this case, $\inf \cJ_N=\infty$. We  claim that 
\bel{w1225e15}
x\big(t;x_0,u(\cd)\big)\neq \xi \qq \forall\, t\in [T_0,\infty)\,, u(\cd)\in \cU_N\,,
\ee
where $\cU_N$ is given in \rf{w1227e1}.

If not, then there would exists $t\geq T_0$ and $u(\cd)\in L^2_\dbF(\dbR^+;\dbR^m)$ such that
\bel{w1225e17}
x\big(t;x_0,u(\cd)\big)=\xi \qq\mbox{ and } \qq \|u(\cd)\|_{ L^2_\dbF(\dbR^+;\dbR^m)}\leq N\,. 
\ee
By \rf{w1225e17}, we know that $u|_{(0,t)}(\cd)$ is an admissible control of Problem $(NP)_{x_0,\xi}^t$, and
$N(t,x_0,\xi)\leq \|u(\cd)\|_{ L^2_\dbF(0,t;\dbR^m)}\leq N$, which means that $t\in \cJ_N$. That leads to a contradiction since in this case
$\cJ_N=\emptyset$. Subsequently, the claim \rf{w1225e15} holds, and 
\beq
T(N,x_0,\xi)=\infty=\inf \cJ_N\,.
\eeq

\no{\bf Case II.} $\cJ_N\neq \emptyset$

For any $t\in \cJ_N$, it follows that $N(t,x_0,\xi)\leq N$. By the unique solvability of $(NP)_{x_0,\xi}^t$, 
we extend the minimal norm control $u^*(\cd)$ to $\dbR^+$ by setting it to zero outside the interval $(0,t)$,
obtaining $\wt u^*(\cd)$.
Then we can easily check that
\bel{w1225e18}
x\big(t;x_0, \wt u^*(\cd)\big)=\xi \q\mbox{and} \q  \|\wt u^*(\cd)\|_{ L^2_\dbF(\dbR^+;\dbR^m)}= \|u^*(\cd)\|_{ L^2_\dbF(0,t;\dbR^m)}=N(t,x_0,\xi)\leq N\,.
\ee
By the fact that $t$ is arbitrary in $\cJ_N$, \rf{w1225e18} and the definition of $T(\cd,x_0,\xi)$, we find that 
\bel{w1225e19}
T(N,x_0,\xi)\leq \inf \cJ_N\,.
\ee

To prove the reversion of \rf{w1225e19}, define
\bel{w1225e20}
\cT_{N,x_0,\xi}\deq \big\{t\geq T_0\,\big|\, \mbox{there exists some } u(\cd)\in \cU_N \mbox{ such that } x\big(t;x_0,u(\cd)\big)=\xi \big\}\,.
\ee
By \rf{w1225e18}, if follows that $\cT_{N,x_0,\xi}\neq \emptyset$. For any $\wt t \in \cT_{N,x_0,\xi}$, there exists a control $\wt u(\cd)$ such that
\bel{w1225e21}
x\big(\wt t; x_0,\wt u(\cd)\big)=\xi \qq\mbox{ and }\qq \|\wt u(\cd)\|_{ L^2_\dbF(\dbR^+;\dbR^m)}\leq N\,,
\ee
by which we can deduce that $\wt u|_{(0,\wt t\,)}(\cd)$ is an admissible control of Problem $(NP)_{x_0,\xi}^{\wt t}$, and 
\beq
N(\wt t ,x_0, \xi)\leq \big\| \wt u|_{(0,\wt t\,)}(\cd)\big\|_{ L^2_\dbF(0,\wt t;\dbR^m)}\leq N\,.
\eeq
That yields that $\wt t \in \cJ_N$ and $\inf \cJ_N\leq \wt t$.
Since $\wt t \in \cT_{N,x_0,\xi} $ is arbitrary, we have that $\inf \cJ_N\leq \inf \cT_{N,x_0,\xi}=T(N,x_0,\xi)$.
That completes the proof.
\end{proof}

In the following, we introduce an auxiliary optimal control problem, denoted as Problem $(\wt{NP})_{x_0,\xi}^T$,
with a minimal value $N_{\tt aux}(T,x_0,\xi)$. Notably, it holds that $N_{\tt aux}(T,x_0,\xi)=N(T, x_0, T)$; see Proposition \ref{w1113p1}. By establishing the continuity of  $N_{\tt aux}(\cd,x_0,\xi)$, we demonstrate the equivalence of concerned  Problems $(TP)_{x_0,\xi}^N$ and $(NP)_{x_0,\xi}^T$.

Under the assumption {\bf (A)} regarding the exact controllability of \rf{sys-1}, and based on \cite{Peng94} 
there exists an invertible matrix $M$ such that 
\beq
DM=\big(I_n, O_{n\times (m-n)}\big)\,.
\eeq
We then set
\bel{w103e1}
BM=\big(K, L\big)\,,\q 
M^{-1} u(\cd)=
\begin{pmatrix}
z(\cd)\\
v(\cd)\\
\end{pmatrix}\,, \q
M^\top M=
\begin{pmatrix}
R& S^\top\\
S & N\\
\end{pmatrix}\,,
\ee
where $K\,,R\in \dbR^{n\times n}$, $L\in \dbR^{n\times (m-n)}$, $N\in \dbR^{(m-n)\times (m-n)}$, $S\in \dbR^{(m-n)\times n}$,
and $z(\cd)\in L^2_\dbF(0,T;\dbR^n)$,  $v(\cd)\in L^2_\dbF(0,T;\dbR^{m-n})$.
By the invertibility of $M$ and  \rf{w103e1}, we find that 
\bel{w1113e1}
N>0\,,\qq R_0\deq R-S^\top N^{-1}S>0\,.
\ee
Thus, we can rewrite \rf{sys-1} as follows 
\bel{sys-02}
\bal
{\mathrm d} x(t)&=\big[Ax(t)+Kz(t)+Lv(t)\big]\rd t+\big[Cx(t)+z(t)\big]\rd W(t)\\
&=\Big\{\big[A-\big(K-LN^{-1}S\big)C\big]x+\big(K-LN^{-1}S\big)\big[Cx(t)+z(t)\big]\\
&\qq +LN^{-1/2}\big[N^{1/2}v(t)+N^{-1/2}Sz(t)\big]\Big\}\rd t+\big[Cx(t)+z(t)\big]\rd W(t)\\
&=: \big[A_0 x(t)+K_0 z_0(t)+L_0 v_0 (t)\big]\rd t+z_0(t) \rd W(t) \qq t\geq 0\,,
\eal
\ee
where
\beq
&A_0=A-\big(K-LN^{-1}S\big)C\,, \q K_0=K-LN^{-1}S\,,\q L_0=LN^{-1/2}\,,\\
\mbox{and}\q & z_0(\cd)=Cx(\cd)+z(\cd)\,,\q v_0(\cd)=N^{1/2}v(\cd)+N^{-1/2}Sz(\cd)\,.
\eeq
Subsequently, by applying \rf{w103e1}, \rf{sys-02} and \rf{w1113e1}, the $L^2$-norm of the control variable $u(\cd)$ in \rf{sys-1}, corresponding to Problem $(NP)_{x_0,\xi}^T$, can be rewritten as
\begin{eqnarray}
&&\|u(\cd)\|^2_{L^2_\dbF(0,T;\dbR^m)}
=\bigg\lan
\begin{pmatrix}
N& S\\
S^\top  & R\\
\end{pmatrix}
\begin{pmatrix}
v(\cd)\\
z(\cd)\\
\end{pmatrix},
\begin{pmatrix}
v(\cd)\\
z(\cd)\\
\end{pmatrix}\bigg\ran_{L^2_\dbF(0,T;\dbR^{m})} \notag\\
&&=\bigg\lan
\begin{pmatrix}
I& 0\\
0  & R-S^\top N^{-1}S\\
\end{pmatrix}
\begin{pmatrix}
N^{1/2}v(\cd)+N^{-1/2}Sz(\cd)\\
Cx(\cd)+z(\cd)-Cx(\cd)\\
\end{pmatrix},
\begin{pmatrix}
N^{1/2}v(\cd)+N^{-1/2}Sz(\cd)\\
Cx(\cd)+z(\cd)-Cx(\cd)\\
\end{pmatrix}\bigg\ran_{L^2_\dbF(0,T;\dbR^{m})} \notag \\
&&=\bigg\lan
\begin{pmatrix}
I& 0\\
0  & R_0\\
\end{pmatrix}
\begin{pmatrix}
v_0(\cd)\\
z_0(\cd)-Cx(\cd)\\
\end{pmatrix},
\begin{pmatrix}
v_0(\cd)\\
z_0(\cd)-Cx(\cd)\\
\end{pmatrix}\bigg\ran_{L^2_\dbF(0,T;\dbR^{m})}\notag \\
&&=\big\lan R_0\big(z_0(\cd)-Cx(\cd)\big), z_0(\cd)-Cx(\cd)\big\ran_{L^2_\dbF(0,T;\dbR^{n})}
 +\| v_0(\cd)\|^2_{L^2_\dbF(0,T;\dbR^{m-n})}\,. \label{w1113e2}
\end{eqnarray}
Now we can propose Problem $(\wt {NP})_{ x_0,\xi}^T$ associated with system \rf{sys-02}, stated as
\beq
N_{\tt aux}(T,x_0,\xi) &\deq \inf\Big\{ \Big[\big\lan R_0\big(z_0(\cd)-Cx(\cd)\big), z_0(\cd)-Cx(\cd)\big\ran_{L^2_\dbF(0,T;\dbR^{n})} +\| v_0(\cd)\|^2_{L^2_\dbF(0,T;\dbR^{m-n})} \Big]^{1/2} \\
&\qq\qq \,\Big|\, \big(z_0(\cd),v_0(\cd)\big)\in L^2_\dbF(0,T;\dbR^n\times \dbR^{m-n}), x\big(T;x_0,z_0(\cd), v_0(\cd)\big)=\xi \Big\}\,.
\eeq
Based on \rf{w1113e2}, it is straightforward to obtain the following connection between $N(T,x_0,\xi)$ and $N_{\tt aux}(T,x_0,\xi)$.
\bp{w1113p1}
Suppose that the assumption {\bf (A)} holds, and $N(T,x_0,\xi)$ resp.~$N_{\tt aux}(T,x_0,\xi)$ are minimal norms of Problems $(NP)_{x_0,\xi}^T$ resp.~$(\wt {NP})_{x_0,\xi}^T$. Then
\bel{w1113e4}
N(T,x_0,\xi)=N_{\tt aux}(T,x_0,\xi)\,.
\ee 
\ep

Below, 
we view \rf{sys-02} with a terminal condition $x(T)=\xi$ as a controlled BSDE, where $v_0(\cd)$ is the control variable and $\big(x(\cd),z_0(\cd)\big)$ 
is the state. We introduce a parameterized LQ problem related to Problem $(\wt {NP})_{x_0,\xi}^T$:\\
{\bf Problem $(BLQ)_\l$.} For any given $\l\in \dbR^n$, search for $v_0^*(\cd)\in L^2(0,T;\dbR^{(m-n)})$ that minimizes the cost functional
\bel{b-cost-functional}
\bal
J_\l\big(v_0(\cd)\big)&=\big\lan R_0\big(z_0(\cd)-Cx(\cd)\big), z_0(\cd)-Cx(\cd)\big\ran_{L^2_\dbF(0,T;\dbR^{n})}\\
&\q+\| v_0(\cd)\|^2_{L^2_\dbF(0,T;\dbR^{m-n})}+2\lan \l, x(0)\ran\,,
\eal
\ee
subject to BSDE
\beq
\lt\{
\bal
&{\mathrm d} x(t)=\big[A_0 x(t)+K_0 z_0(t)+L_0 v_0 (t)\big]\rd t+z_0(t) \rd W(t) \qq t\in [0,T]\,,\\
&x(T)=\xi\,.
\eal
\rt.
\eeq

Differing from Problem $(\wt{NP})_{x_0,\xi}^T$, Problem $(BLQ)_\l$ is a standard LQ problem without any constraints.
The following discussion focuses on the unique solvability of Problem $(BLQ)_\l$, and its proof follows standard argument; see e.g., \cite[Theorem 3.3]{Lim-Zhou01}, \cite[Theorem 4.2]{Bi-Sun-Xiong20}.
\bp{w1229p1}
For any $\l\in \dbR^n$, the unique solvability of Problem $(BLQ)_\l$ 
is equivalent to the unique solvability of the following coupled forward-backward SDE
\bel{fbsde}
\lt\{
\bal
&{\mathrm d} x^*(t)=\big[A_0 x^*(t)+K_0 z_0^*(t)+L_0 v_0^* (t)\big]\rd t+z_0^*(t) \rd W(t) \qq t\in [0,T]\,,\\
&{\mathrm d} y^*(t)=\big[-A_0^\top y^*(t)+C^\top R_0 C x^*(t) - C^\top R_0 z_0^*(t)\big]\rd t\\
&\qq\qq+ \big[-K_0^\top y^*(t)-R_0C x^*(t) + R_0 z_0^*(t) \big]  \rd W(t) \qq t\in [0,T]\,,\\
&x(T)=\xi\,, y^*(0)=\l\,,
\eal
\rt.
\ee
with the optimality condition
\bel{optimality-condition}
v_0^*(\cd)=L_0^\top y^*(\cd)\,.
\ee
\ep

In LQ theory, Riccati equations play an important role in feedback representations of the optimal controls; 
see e.g., \cite{Yong-Zhou99, Lim-Zhou01,Lv19}.
In the below, we also attempt to derive the optimal control's representation of Problem $(\wt{NP})_{x_0,\xi}^T$. Hence,
 we introduce a Riccati equation
\bel{Riccati-1}
\lt\{
\bal
& \dot \Si(t)-A_0\Si(t)-\Si(t)A_0^\top +L_0L_0^\top-\Si(t)C^\top R_0 C \Si(t)\\
&\q+ \big(K_0-\Si(t)C^\top R_0\big) \big[I+\Si(t)R_0\big]^{-1}\Si(t)\big(K_0-\Si(t)C^\top R_0\big)^\top=0 \qq t\in [0,T]\,,\\
&\Si(T)=0\,.
\eal
\rt.
\ee
By following the proof in \cite{Lim-Zhou01, Li-Sun-Xiong19}, one can deduce that Riccati equation \rf{Riccati-1} has a unique 
solution $\Si(\cd)\in C([0,T];\dbS^n_+)$, where $\dbS^{n}_+$ is the set of $n\times n$ positive semidefinite matrices. 
By 
setting 
\bel{w1229e5}
\bal
&\h K(\cd)=\big(K_0-\Si(\cd)C^\top R_0\big) \big[I+\Si(\cd)R_0\big]^{-1}\,,\\
\mbox{and }\q&\h A(\cd)=A_0+\Si(\cd)C^\top R_0 C+\h K(\cd )\Si(\cd)R_0 C\,,
\eal
\ee
we introduce a BSDE 
\bel{bsde-3}
\lt\{
\bal
& {\mathrm d} \f(t)=\big[\h A(t) \f(t)+\h K(t)\b(t)\big]\rd t +\b(t) \rd W(t) \qq t\in [0,T]\,,\\
&\f(T)=-\xi\,,
\eal
\rt.
\ee
which has a unique solution $\big(\f(\cd),\b(\cd)\big)\in L^2_\dbF\big(\Omega;C([0,T];\dbR^n)\big)\times L^2_\dbF(0,T;\dbR^n)$. 

The following result concerning BSDE \rf{bsde-3}, derived from the uniqueness of its solution, will be
applied in Theorem \ref{w1230t1}.
\bl{w1123l1}
Suppose that $\xi\in L^2_{\mf_{T_0}}(\O;\dbR^n)$ with $T_0<T$, and let $\big(\f(\cd), \b(\cd)\big)$ solve BSDE \rf{bsde-3}. Then, it holds that
\beq
\big(\f(t),\b(t)\big)=\big(-\Pi_0^{-1}(t)\Pi_0(T)\xi, 0\big)\qq \forall\, t\in [T_0, T]\,,
\eeq
where $\Pi_0(\cd)$ is the solution of the following ODE
\beq
\lt\{
\bal
& \dot \Pi_0(t)=-\Pi_0(t)\h A(t) \qq t\in [0,T]\,,\\
&\Pi_0(0)=I_n\,.
\eal
\rt.
\eeq
\el

Since system \rf{sys-1} is exactly controllable, we know that 
system \rf{sys-02}  is also exactly controllable. 
To derive the invertibility of $\Si(0)$, which is pivotal for constructing the optimal tuple of Problem $(\wt{NP})_{x_0,\xi}^T$,
we introduce another system:
\bel{sys-05}
\bal
 {\mathrm d} \h x(t)=\big[\h A(t) \h x(t)+\h K(t)\h z(t)+\h L(t) \h v(t)\big]\rd t +\h z(t) \rd W(t) \qq t\in [0,T]\,,
\eal
\ee
where 
\bel{w1229e6}
\bal
\h L(\cd)&=\big(\h L_1(\cd)\,,\h L_2(\cd)\big)\\
&\deq
\begin{pmatrix}
L_0\,, & K_0\big[I+\Si(\cd)R_0\big]^{-1}\Si(\cd)R_0^{1/2}+\Si(\cd)C^\top \big[I+R_0\Si(\cd)\big]^{-1}R_0^{1/2}   
\end{pmatrix}\,,
\eal
\ee
and
\bel{w1229e7}
\h z(\cd)=z_0(\cd)\,,\q
\h v(\cd)=
\begin{pmatrix}
v_0(\cd)\\
R_0^{1/2}z_0(\cd)-R_0^{1/2}Cx(\cd)
\end{pmatrix}\,.
\ee

\bp{w1229p2}
Under the assumption {\rm \bf (A)},
system \rf{sys-05} is exactly controllable. Moreover, the Gramian matrix $W(T)$ of \rf{sys-05} is invertible
and
\bel{w1229e8}
W(T)=\Si(0)\,.
\ee
Here 
\bel{w1229e9}
W(T)=\me\int_0^T \Pi(t)\h L (t)\h L(t)^\top \Pi(t)^\top \rd t \,,
\ee
$\Pi(\cd)$ solves
\bel{Pi-eq}
\lt\{
\bal
& {\mathrm d} \Pi(t)=-\Pi(t)\h A(t)\rd t-\Pi(t) \h K(t) \rd W(t) \qq t\in [0,T]\,,\\
&\Pi(0)=I_n\,,
\eal
\rt.
\ee
and
$\Si(\cd)$ is the solution to Riccati equation \rf{Riccati-1}.
\ep

\begin{proof}
{\bf (1)}
Under the assumption {\bf (A)},
for any initial state $x_0\in \dbR^n$ and any target $\xi\in L^2_{\mf_T}(\Omega;\dbR^n)$, controllability of system \rf{sys-1}, or equivalently \rf{sys-02}, guarantees the existence of $\big(z_0(\cd), v_0(\cd)\big)$ such that
\beq
x(0)=x_0\,,\q x(T)=\xi\,.
\eeq 

By applying \rf{w1229e5}--\rf{w1229e7}, system \rf{sys-05} with $\h x(0)=x_0$ turns to
\bel{sys-051}
\lt\{
\bal
& {\mathrm d} \h x(t)=\Big[A_0 \h x(t)+ K_0z_0(t)+ L_0  v_0(t)+\h L_2(t)R_0^{1/2}C\big(\h x(t)-x(t)\big)\Big]\rd t\\
&\qq\qq + z_0(t) \rd W(t) \qq t\in [0,T]\,,\\
&\h x(0)=x_0\,.
\eal
\rt.
\ee
Thus, unique solvability and linearity of \rf{sys-02} with initial condition $x(0)=x_0$, along with \rf{sys-051}, yield that
\beq
\h x(\cd)=x(\cd)\,, \q\mbox { and } \q \h x(T)=x(T)=\xi\,.
\eeq
That settles the exact controllability of system \rf{sys-05}.

\ms
{\bf (2)} In this step, we prove that the Gramian matrix $W(T)$ 
is invertible and $W(T)=\Si(0)$.
The invertibility of $W(T)$ can be derived in the same vein as that in \cite[Theorem 2]{Liu-Peng10}.
In the below, we only derive \rf{w1229e8}. To do that, It\^o's formula to $\Pi(t)\Si(t)\Pi(t)^\top$ yields that
\begin{eqnarray}
&&\rd (\Pi(t)\Si(t)\Pi(t)^\top) \notag \\
&&=-\Pi(t)\bigg\{\Big[A_0+K_0\big[I+\Si R_0\big]^{-1}\Si R_0C+\Si C^\top R_0 \big[I+\Si R_0\big]^{-1} C\Big] \Si \notag \\
&&\q+\Big[-A_0\Si-\Si A_0^\top +L_0L_0^\top-\Si C^\top R_0 C \Si \notag \\
&&\qq+ \big(K_0-\Si C^\top R_0\big) \big[I+\Si R_0\big]^{-1}\Si \big(K_0-\Si C^\top R_0\big)^\top\Big] \notag \\
&&\q+\Si \Big[A_0+K_0\big[I+\Si R_0\big]^{-1}\Si R_0C+\Si C^\top R_0 \big[I+\Si R_0\big]^{-1} C\Big]^\top  \notag \\
&&\q- \big(K_0-\Si C^\top R_0\big) \big[I+\Si R_0\big]^{-1}\Si  \big[I+R_0\Si \big]^{-1} \big(K_0-\Si C^\top R_0\big) ^\top\bigg\}\Pi(t)^\top \rd t  \notag \\
&&\q+\cds \rd W(t) \notag \\
&&=-\Pi(t)\bigg\{ \Si C^\top R_0 \big[I+\Si R_0\big]^{-1} C \Si + L_0L_0^\top \notag \\
&&\q+ K_0 \big[I+\Si R_0\big]^{-1}\Si K_0^\top -\Si C^\top R_0 \big[I+\Si R_0\big]^{-1}\Si K_0^\top \notag \\
&&\q+\Si C^\top R_0\Si \big[I+\Si R_0\big]^{-\top}K_0^\top \notag  \\
&&\q -K_0 \big[I+\Si R_0\big]^{-1}\Si  \big[I+R_0\Si \big]^{-1} K_0^\top
+K_0 \big[I+\Si R_0\big]^{-1}\Si  \big[I+R_0\Si \big]^{-1}R_0C\Si \notag  \\
&&\q + \Si C^\top R_0\big[I+\Si R_0\big]^{-1}\Si  \big[I+R_0\Si \big]^{-1} K_0^\top \notag \\
&&\q- \Si C^\top R_0\big[I+\Si R_0\big]^{-1}\Si  \big[I+R_0\Si \big]^{-1}R_0C\Si 
\bigg\}\Pi(t)^\top \rd t \notag  \\
&&\q+\cds \rd W(t) \notag \\
&&=:-\Pi(t)\bigg\{\sum_{i=1}^9I_i(t) \bigg\}\Pi(t)^\top \rd t +\cds \rd W(t)\,. \label{w1229e10}
\end{eqnarray}
We now proceed to calculate the summation. Based on the following calculation
\beq
I_3+I_6&=K_0 \big[I+\Si R_0\big]^{-1}\Si R_0\Si  \big[I+R_0\Si \big]^{-1} K_0^\top\,,
\eeq
\beq
I_1+I_9&= \Si C^\top R_0\big[I+\Si R_0\big]^{-1}\big[I+\Si R_0 \big]^{-1}C\Si\\
&= \Si C^\top \big[I+R_0\Si\big]^{-1}R_0\big[I+\Si R_0 \big]^{-1}C\Si\,,
\eeq
and
\beq
(I_4+I_8)+I_5&=-\Si C^\top R_0\big[I+\Si R_0\big]^{-1}\Si R_0\Si \big[I+R_0\Si \big]^{-1} K_0^\top+I_5\\
&=\Si C^\top R_0\big[I+\Si R_0\big]^{-1}\Si \big[I+R_0\Si \big]^{-1} K_0^\top\\
&=\Si C^\top \big[I+ R_0\Si\big]^{-1}R_0 \Si \big[I+R_0\Si \big]^{-1} K_0^\top\,,
\eeq
we find that
\bel{w1229e11}
\sum_{i=1}^9I_i =L_0L_0^\top+\h L_2 \h L_2^\top=\h L \h L^\top\,.
\ee
By combining with \rf{w1229e10}, \rf{w1229e11} and \rf{w1229e9}, we conclude that
\beq
\Si(0)=\me\int_0^T\Pi(t)\sum_{i=1}^9I_i(t) \Pi(t)^\top\rd t=\me\int_0^T\Pi(t) \h L(t) \h L(t)^\top \Pi(t)^\top\rd t=W(T)\,.
\eeq
That completes the proof.
\end{proof}

Since $\Si(0)$ is invertible, we can specify a parameter
\bel{w1229e12}
\l^*=-\Si(0)^{-1}\big(x_0-\me\big[\Pi(T)\xi\big]\big)\,,
\ee
and then we can derive the unique optimal tuple $\big(x_{\l^*}^*(\cd),z_{\l^*}^*(\cd),v_{\l^*}^*(\cd)\big)$ of Problem $(BLQ)_{\l^*}$.
The following result shows that $\big(z_{\l^*}^*(\cd),v_{\l^*}^*(\cd)\big)$ is the optimal control of Problem $(\wt{NP})_{x_0,\xi}^T$.

\bp{w1229p4}
Suppose that $\l^*$ is given by \rf{w1229e12}. Then
for  the unique optimal tuple $\big(x_{\l^*}^*(\cd),z_{\l^*}^*(\cd),v_{\l^*}^*(\cd)\big)$ of Problem $(BLQ)_{\l^*}$, 
it holds that
\bel{w1229e13}
\lt\{
\bal
&x_{\l^*}^*(\cd)=-\Si(\cd)y^*_{\l^*}(\cd)-\f(\cd)\,,\\
&z_{\l^*}^*(\cd)=\big[I+\Si(\cd) R_0\big]^{-1}\Big[\Si(\cd)\big(K_0^\top-R_0C\Si(\cd) \big)y^*_{\l^*}(\cd)-\Si(\cd) R_0C\f(\cd)-\b(\cd)\Big]\,,\\
&v_{\l^*}^*(\cd)= L_0^\top y^*_{\l^*}(\cd)\,,
\eal
\rt.
\ee
where $y^*_{\l^*}(\cd)$ solves the following SDE
\bel{sde-1}
\lt\{
\bal
& {\mathrm d} y(t)=\Big[-\h A^\top(t)  y(t)- C^\top R_0 \big[I+\Si(t)R_0 \big]^{-1}\big(C\f(t)-\b(t)\big)\Big]\rd t\\
&\qq\q +\Big[-\h K(t)^\top  y(t)+ R_0 \big[I+\Si(t)R_0 \big]^{-1}\big(C\f(t)-\b(t)\big)\Big]\rd W(t) \qq t\in [0,T]\,,\\
&y (0)=\l^*=-\Si(0)^{-1}\big(x_0-\me\big[\Pi(T)\xi\big]\big)\,.
\eal
\rt.
\ee
Furthermore, $\big(x_{\l^*}^*(\cd), z_{\l^*}^*(\cd),v_{\l^*}^*(\cd)\big)$ minimizes $N_{\tt aux}(T,x_0,\xi)$ associated to Problem $(\wt {NP})_{x_0,\xi}^T$.
\ep

\begin{proof}
It is straightforward to verify that the quadruple $\big(x^*_{\l^*}(\cd),z^*_{\l^*}(\cd),v^*_{\l^*}(\cd),y^*_{\l^*}(\cd)\big)$
satisfies \rf{fbsde} and \rf{optimality-condition}. Hence, by Proposition \ref{w1229p1}, it uniquely solves Problem $(BLQ)_{\l^*}$.

By applying It\^o's formula to $\Pi(\cd)\f(\cd)$, we find that
\beq
\f(0)= \Pi(0)\f(0)= \me\big[\Pi(T)\f(T)\big]= -\me\big[\Pi(T)\xi\big]\,,
\eeq
which, together with \rf{w1229e13}, leads to
\beq
x^*_{\l^*}(0)=x_0-\me\big[\Pi(T)\xi\big]-\f(0)=x_0\,.
\eeq
Therefore, $\big(x_{\l^*}^*(\cd),z_{\l^*}^*(\cd),v_{\l^*}^*(\cd)\big)$ is an admissible tuple of Problem $(\wt {NP})_{x_0,\xi}^T$. 
Now, consider another admissible tuple $\big(x(\cd),z(\cd),v(\cd)\big)$ of Problem $(\wt {NP})_{x_0,\xi}^T$.
Clearly, $v(\cd)$ is an admissible control for Problem $(BLQ)_{\l^*}$. By the optimality of $\big(x_{\l^*}^*(\cd),z_{\l^*}^*(\cd),v_{\l^*}^*(\cd)\big)$, we have
\beq
J_{\l^*}\big(v^*_{\l^*}(\cd)\big)
&=\big\lan R_0\big(z^*_{\l^*}(\cd)-Cx^*_{\l^*}(\cd)\big), z^*_{\l^*}(\cd)-Cx^*_{\l^*}(\cd)\big\ran_{L^2_\dbF(0,T;\dbR^{n})}
+\| v^*_{\l^*}(\cd)\|^2_{L^2_\dbF(0,T;\dbR^{m-n})}+2\lan \l^*, x_0\ran\\
&\leq \big\lan R_0\big(z(\cd)-Cx(\cd)\big), z(\cd)-Cx(\cd)\big\ran_{L^2_\dbF(0,T;\dbR^{n})}
+\| v(\cd)\|^2_{L^2_\dbF(0,T;\dbR^{m-n})}+2\lan \l^*, x_0\ran\\
&=J_{\l^*}\big(v(\cd)\big)\,,
\eeq
which implies that 
\beq
&\big\lan R_0\big(z^*_{\l^*}(\cd)-Cx^*_{\l^*}(\cd)\big), z^*_{\l^*}(\cd)-Cx^*_{\l^*}(\cd)\big\ran_{L^2_\dbF(0,T;\dbR^{n})}
+\| v^*_{\l^*}(\cd)\|^2_{L^2_\dbF(0,T;\dbR^{m-n})}\\
&\qq\leq \big\lan R_0\big(z(\cd)-Cx(\cd)\big), z(\cd)-Cx(\cd)\big\ran_{L^2_\dbF(0,T;\dbR^{n})}
+\| v(\cd)\|^2_{L^2_\dbF(0,T;\dbR^{m-n})}\,.
\eeq
Thus, $\big(x_{\l^*}^*(\cd),z_{\l^*}^*(\cd),v_{\l^*}^*(\cd)\big)$ is an optimal tuple of Problem $(\wt {NP})_{x_0,\xi}^T$. 
That completes the proof.
\end{proof}

\br{w1229r1}
Proposition \ref{w1229p4} provides a scheme to derive the minimal norm control of Problem $(NP)_{x_0,\xi}^T$:
\begin{enumerate}

\item [{\rm (1)}] Solve $\Si(\cd)$, $\big(\f(\cd),\b(\cd)\big)$  resp.~$\Pi(\cd)$ by \rf{Riccati-1}, \rf{bsde-3} resp.~\rf{Pi-eq};

\item [{\rm (2)}] Derive $\l^*$ by \rf{w1229e12}, and then solve SDE \rf{sde-1};

\item [{\rm (3)}] Derive $\big(x_{\l^*}^*(\cd),z_{\l^*}^*(\cd),v_{\l^*}^*(\cd)\big)$ by \rf{w1229e13};

\item [{\rm (4)}] Derive the minimal norm control of Problem $(NP)_{x_0,\xi}^T$ by invertible transformation \rf{sys-02}
and \rf{w103e1}.

\end{enumerate}
\er

The following result is on the continuity of the minimal norm function $N(\cd,x_0,\xi)$, which is crucial in deriving the 
equivalence between minimal time and norm problems.

\bt{w1230t1}
Under the assumption {\rm \bf (A)},
for any given $x_0\in \dbR^n$ and $\xi \in L^2_{\mf_{T_0}}(\Omega;\dbR^n)$ satisfying \rf{w1226e1},
$N^2(\cd,x_0, \xi)$,  the square of minimal norm function,  is locally $\frac 1 4$-H\"older continuous on $[T_0, \infty)$. 
\et

\begin{proof}

Based on the connection between $N(T,x_0,\xi)$ and $N_{\tt aux}(T,x_0,\xi)$ stated in Proposition \ref{w1113p1}, it suffices to establish the locally 
$\frac 1 4$-H\"older continuity for $N_{\tt aux}^2(\cd,x_0,\xi)$. 

In order to simplify notations,
for any $T_1,T_2\in [T_0, \infty)$ and $T_1\leq T_2$, we will write the optimal tuple of Problems $(\wt{NP})_{x_0,\xi}^{T_1}$
resp.~$(\wt{NP})_{x_0,\xi}^{T_2}$ by $\big(x_1(\cd),z_1(\cd),v_1(\cd)\big)$ resp.~$\big(x_2(\cd),z_2(\cd),v_2(\cd)\big)$.
Similarly we can introduce $\Si_i(\cd)\,,\h A_i (\cd)\,,\h K_i(\cd)$, $\big(\f_i(\cd),\b_i(\cd)\big)$, $\Pi_i(\cd)$,  $y_i(\cd)$ for $i=1\,,2$\,. Denote by $\d\deq T_2-T_1$, and  for $i=1\,,2$
by $f_i(\cd)\deq R_0 \big[I+\Si_i(\cd)R_0 \big]^{-1}\big(C\f_i(\cd)-  \b_i(\cd)\big)$.

\ms
{\bf (1)} Property of $\Si_i(\cd)$ to Riccati equation \rf{Riccati-1}

It is obvious that $\Si_1(t)=\Si_2(t+\d)$ for any $t\in [0,T_1]$. Since $\big[I+\Si_i(\cd)R_0 \big]^{-1}$ is uniformly bounded on $[0,T_i]$, and then for any $t\in [0,T_1]$
\bel{w1230e1}
\|\Si_1(t)-\Si_2(t)\|
=\Big\|\int_t^{t+\d} \dot \Si_2(s) \rd s\Big\|\leq \cC(T_2)\d\,,
\ee
and
\bel{w1230e2}
\bal
&\big\|\big[I+\Si_1(t)R_0 \big]^{-1}-\big[I+\Si_2(t)R_0 \big]^{-1}\|\\
&\q=\big\|\big[I+\Si_1(t)R_0 \big]^{-1}\big\| \times\big\|\big(\Si_1(t)-\Si_2(t)\big)R_0\big\| \times \big\|\big[I+\Si_2(t)R_0 \big]^{-1}\big\|\\
&\q\leq \cC(T_1,T_2)\d\,.
\eal
\ee
Besides, Proposition \ref{w1229p2} yields that
\bel{w1230e4}
0<\cC(T_i)I\leq \Si_i(0) \qq i=1\,,2\,.
\ee

{\bf (2)} Property of $\big(\f_i(\cd),\b_i(\cd)\big)$,  $\Pi_i(\cd)$ and $f_i(\cd)$

Relying on the property of $\Si_i(\cd)$ derived in Step (1),
we can deduce that 
\bel{w103e4}
\bal
&\sup_{t\in [0,T_i]}\big[\big\|\h A_i(t)\big\|+\big\|\h K_i(t)\big\|+\me\big(\|\Pi_i(t)\|^2\big)\big]\leq \cC(T_i) \qq i=1\,,2\,,\\
&\sup_{t\in [0,T_1]}\big[\big\|\h A_1(t)-\h A_2(t)\big\|+\big\|\h K_1(t)-\h K_2(t)\big\|\big]\leq \cC(T_1,T_2)\d\,,
\eal
\ee
and 
then 
\bel{w1230e5}
\me\Big[\sup_{t\in[0,T_i]}\|\f_i(t)\|^2+\int_0^{T_i}\|\b_i(t)\|^2\rd t \Big]\leq \cC(T_i) \qq i=1\,,2\,.
\ee
By \rf{w1230e5} and Lemma \ref{w1123l1}, we have
\bel{w1123e1}
\me\int_{T_1}^{T_2}\|f_2(t)\|^2\rd t
=\me\int_{T_1}^{T_2} \big\| R_0 \big[I+\Si_2(t)R_0 \big]^{-1}C\f_2(t)\big\|^2\rd t
\leq \cC(T_2)\d\,.
\ee
Applying \rf{w1230e2} and BSDE \rf{bsde-3}, it follows that
\bel{w1230e6}
\bal
&\me\int_0^{T_1}\big(\|\f_1(t)-\f_2(t)\|^2+\|\b_1(t)-\b_2(t)\|^2\big)\rd t \\
&\q\leq \cC(T_1,T_2)\me\Big\{\|-\xi-\f_2(T_1)\|^2\\
&\qq\q+\int_0^{T_1}\big[\big\|\big(\h A_1(t)-\h A_2(t)\big)\f_2(t)\big\|^2+\big\|\big(\h K_1(t)-\h K_2(t)\big)\b_2(t)\big\|^2\big] \rd t\Big\}\\
&\q\leq \cC(T_1,T_2)\d\,,
\eal
\ee
and then
\bel{w105e1}
\me\int_0^{T_1}\|f_1(t)-f_2(t)\|^2 \rd t \leq \cC(T_1,T_2)\d\,.
\ee
Note that to estimate $\me \big(\|-\xi-\f_2(T_1)\|^2\big)$ in \rf{w1230e6}, we apply the fact that $-\xi-\f_2(T_1)$ is $\mf_{T_1}$-measurable and 
\beq
-\xi-\f_2(T_1)=\me\big(-\xi-\f_2(T_1)\,\big|\,\mf_{T_1}\big)
=\me\Big(\int_{T_1}^{T_2} \h A_2(t)\f_2(t)+\h K_2(t)\b_2(t)\rd t \,\Big|\,\mf_{T_1}\Big)\,.
\eeq
Similarly, \rf{w103e4} and SDE \rf{Pi-eq} lead to
\bel{w103e5}
\me\big[\|\Pi_2(T_1)-\Pi_2(T_2)\|^2\big]\leq \cC(T_1,T_2)\d\,,
\ee
and
\bel{w103e6}
\bal
&\me\big[\|\Pi_1(T_1)-\Pi_2(T_1)\|^2\big]\\
&\q\leq \cC(T_1,T_2)\me \int_0^{T_1}\big[\big\|\Pi_2(t)\big(\h A_1(t)-\h A_2(t)\big)\big\|^2+\big\|\Pi_2(t)\big(\h K_1(t)-\h K_2(t)\big)\big\|^2\big] \rd t\\
&\q\leq \cC(T_1,T_2)\d^2\,.
\eal
\ee

{\bf (3)} Difference between  $\|y_1(\cd)\|^2_{L^2_\dbF(0,T_1;\dbR^n)}$ and $\|y_2(\cd)\|^2_{L^2_\dbF(0,T_2;\dbR^n)}$

Regularity of solution to SDE \rf{sde-1} with $T=T_1\,, T_2$ and \rf{w1230e4} lead to
\bel{w1230e7}
\bal
&\me\Big[\sup_{t\in [0,T_i]}\|y_i(t)\|^2\Big]\\
&\q\leq \cC(T_i)\Big[\big\|\Si_i(0)^{-1}\big(x_0-\me\big[\Pi_i(T_i)\xi\big]\big)\big\|^2
+\me\int_0^{T_i}\|f_i(t)\|^2\rd t\Big]\\
&\q\leq \cC(T_i) \qq i=1\,,2\,;
\eal
\ee
and \rf{w1230e1}--\rf{w103e4}, \rf{w105e1}--\rf{w103e6} yield that
\bel{w1230e8}
\bal
&\me\int_0^{T_1}\|y_1(t)-y_2(t)\|^2\rd t\\
&\leq \cC(T_1,T_2)\Big\{\big\|\Si_1(0)^{-1}-\Si_2(0)^{-1}\big\|^2
+\big[\me\big(\|\Pi_1(T_1)-\Pi_2(T_2)\|^2\big) \me\big(\|\xi\|^2\big)\big]^{1/2}\\
&\qq+\me\int_0^{T_1}\big[\big\|\big(\h A_1(t)-\h A_2(t)\big)^\top y_2(t)\big\|^2
+\big\|\big(\h K_1(t)-\h K_2(t)\big)^\top y_2(t)\big\|^2
+\|f_1(t)-f_2(t)\|^2\big]\rd t\Big\}\\
&\leq \cC(T_1,T_2)\d^{1/2} \,.
\eal
\ee
By combining with estimates \rf{w1230e7} and \rf{w1230e8}, we arrive at
\bel{w1230e9}
\bal
&\big|\|y_2(\cd)\|^2_{L^2_\dbF(0,T_2;\dbR^n)}-\|y_1(\cd)\|^2_{L^2_\dbF(0,T_1;\dbR^n)}\big|\\
&=\Big|\me\int_0^{T_1}\lan y_2(t)-y_1(t), y_2(t)+y_1(t)\ran \rd t+\me\int_{T_1}^{T_2}\|y_2(t)\|^2\rd t\Big|\\
&\leq\cC\Big\{ \me\Big[ \sup_{t\in [0,T_1]}\big(\|y_1(t)\|^2+\|y_2(t)\|^2\big)\Big] \me\int_0^{T_1}\|y_1(t)-y_2(t)\|^2\rd t\Big\}^{1/2}
 +\me\int_{T_1}^{T_2}\|y_2(t)\|^2\rd t \\
&\leq \cC(T_1,T_2)\d^{1/4} \,.
\eal
\ee

\ms
{\bf (4)} Difference between  $\|v_1(\cd)\|^2_{L^2_\dbF(0,T_1;\dbR^{m-n})}$ and $\|v_2(\cd)\|^2_{L^2_\dbF(0,T_2;\dbR^{m-n})}$

Optimality condition $v_i(\cd)=L_0^\top y_i(\cd)$ and \rf{w1230e9} imply that
\beq
\big|\|v_2(\cd)\|^2_{L^2_\dbF(0,T_2;\dbR^{m-n})}-\|v_1(\cd)\|^2_{L^2_\dbF(0,T_1;\dbR^{m-n})}\big| \leq \cC(T_1,T_2)\d^{1/4}\,.
\eeq

\ms
{\bf (5)} Difference between  $\|z_1(t)-Cx_1(\cd)\|^2_{L^2_\dbF(0,T_1;\dbR^n)}$ and $\|z_2(\cd)-Cx_2(\cd)\|^2_{L^2_\dbF(0,T_2;\dbR^n)}$

Relying on \rf{w1229e13}, we can derive 
\beq
z_i(t)-Cx_i(t)=\big[\Si_i(t)\h K_i^\top(t)+C\Si_i(t)\big] y_i(t)+R_0^{-1}f_i(t)\q\ae\, t\in(0,T_i)\,.
\eeq
Then, following the same vein as in Step (3) and using estimates in the first three steps, we can conclude that
\beq
\bal
&\big|\|z_2(\cd)-Cx_2(\cd)\|^2_{L^2_\dbF(0,T_2;\dbR^n)}-\|z_1(\cd)-Cx_1(\cd)\|^2_{L^2_\dbF(0,T_1;\dbR^n)}\big|\\
&\leq \cC\bigg\{ \me\sup_{t\in [0,T_1]}\Big[\sum_{i=1}^2\|\Si_i(t)\|^2\big(\|\h K_i(t)\|+\|C\|\big)^2\|y_i(t)\|^2\Big] 
+\sum_{i=1}^2\me\int_0^{T_1}\|f_i(t)\|^2\rd t\bigg\}^{1/2}\\
&\qq\times \bigg\{\me\int_0^{T_1}\|\Si_1(t)-\Si_2(t)\|^2\sup_{t\in[0,T_1]}\big[\big(\|\h K_1(t)\|^2+\|C\|^2\big)\|y_1(t)\|^2\big]\\
&\qq\qq+\|\h K_1(t)-\h K_2(t)\|^2\sup_{t\in[0,T_1]}\big[\|\Si_2(t)\|^2\|y_1(t)\|^2\big]\\
&\qq\qq+\|y_1(t)-y_2(t)\|^2\sup_{t\in[0,T_1]}\big[\|\Si_2(t)\|^2\big(\|\h K_1(t)\|^2+\|C\|^2\big)\big]\\
&\qq\qq+\|f_1(t)-f_2(t)\|^2\rd t\bigg\}^{1/2}\\
&\q +\cC\, \me\int_{T_1}^{T_2}\sup_{t\in[0,T]}\Big[\|\Si_2(t)\|^2\big(\|\h K_2(t)\|+\|C\|\big)^2\|y_2(t)\|^2\Big]+\|f_2(t)\|^2\rd t \\
&\leq \cC(T_1,T_2)\d^{1/4} \,.
\eal
\eeq

\ms
Finally, combining the above two estimates in Steps (4) and (5), we deduce that
\beq
&\big|N_{\tt aux}^2(T_2,x_0,\xi)-N_{\tt aux}^2(T_1,x_0,\xi)\big|\\
&\qq\leq \cC\Big[\big|\|z_2(\cd)-Cx_2(\cd)\|^2_{L^2_\dbF(0,T_2;\dbR^n)}-\|z_1(\cd)-Cx_1(\cd)\|^2_{L^2_\dbF(0,T_1;\dbR^n)}\big|\\
&\qq\qq+\big|\|v_2(\cd)\|^2_{L^2_\dbF(0,T_2;\dbR^{m-n})}-\|v_1(\cd)\|^2_{L^2_\dbF(0,T_1;\dbR^{m-n})}\big|\Big]\\
&\qq\leq \cC(T_1,T_2)\d^{1/4} \,.
\eeq
That completes the proof.
\end{proof}

The following result builds a relationship between $N(T,x_0,\xi)$ and $T(N,x_0,\xi)$.
\bp{w1225p1}
Under the assumption {\rm \bf (A)} and \rf{w1226e1},
 for each $N\geq 0$ with $T_0<T(N,x_0,\xi)<\infty$, it holds that
\bel{w1225e22}
N\big(T(N,x_0,\xi),x_0,\xi\big)=N\,.
\ee
\ep

\begin{proof}
By Lemma \ref{w1225l3}, it follows that $\cJ_N\neq \emptyset$. Then by the definition of $\cJ_\cd$ and $T(\cd,x_0,\xi)$, 
there exists a sequence $\{t_k\}_{k=1}^\infty\subset [T_0,\infty)$ such that
\beq
\lim_{k\to \infty}t_k=T(N,x_0,\xi) \qq\mbox{ and } \qq N(t_k,x_0,\xi)\leq N \qq\forall\, k\in \dbN^+\,.
\eeq
On the other hand, continuity of $N(\cd,x_0,\xi)$ derived in Theorem \ref{w1230t1} implies that
\bel{w1225e23}
N\big(T(N,x_0,\xi),x_0,\xi\big)\leq N\,.
\ee

We proceed to establish the reverse of \rf{w1225e23}. By contradiction, we suppose that it was not true. Then, by \rf{w1225e23}, we would
have $N\big(T(N,x_0,\xi),x_0,\xi\big)< N$. 
Given the continuity of $N(\cd,x_0,\xi)$ at $T(N,x_0,\xi)$, 
there exists a constant $\d_0\in \big(0, T(N,x_0,\xi)-T_0\big)$ such that $N\big(T(N,x_0,\xi)-\d_0,x_0,\xi\big)< N$. 
By the definition of $\cJ_N$, this implies that $T(N,x_0,\xi)-\d_0\in \cJ_N$, which contradicts Lemma \ref{w1225l3}. This completes the proof of assertion \rf{w1225e22}.
\end{proof}

By applying Proposition \ref{w1225p1}, we can settle the equivalence of concerned Problems $(TP)_{x_0,\xi}^N$ and $(NP)_{x_0,\xi}^T$.

\bt{w1225t2}
Let the assumption {\rm \bf (A)} and \rf{w1226e1} hold, and
suppose that $$(N,T)\in \big\{\big(N, T(N,x_0,\xi)\big)\,\big|\, T(N,x_0,\xi)\in(T_0,\infty)\big\}\,.$$ Then Problem $(TP)_{x_0,\xi}^N$
is equivalent to Problem $(NP)_{x_0,\xi}^T$.
\et

\begin{proof}
We divide the proof into three steps.

{\bf (1)}  We claim that: Problems $(NP)^T_{x_0,\xi}$ and $(TP)^N_{x_0,\xi}$ are solvable. 

By Lemma \ref{w727l1}, there exists a unique minimal norm control $u^*(\cd)\in L^2_\dbF(0,T;\dbR^m)$ to Problem
$(NP)_{x_0,\xi}^T$ and
\beq
\|u^*(\cd)\|_{L^2_\dbF(0,T;\dbR^m)}=N(T,x_0,\xi)\,, \qq x\big(T;x_0,u^*(\cd)\big)=\xi\,.
\eeq
These two statements lead to the solvability of Problem $(NP)^{T}_{x_0,\xi}$.
Extend $u^*(\cd)$ to $[0,\infty)$ by setting it to be zero on $(T,\infty)$, and denote the extension by $\wt u^*(\cd)$.
It is obvious that 
\beq
&\qq\|\wt u^*(\cd)\|_{L^2_\dbF(\dbR^+;\dbR^m)}=\|u^*(\cd)\|_{L^2_\dbF(0,T;\dbR^m)}=N(T,x_0,\xi) \\
\mbox{and} &\qq
x\big(T(N,x_0,\xi); x_0,\wt u^*(\cd)\big)=x\big(T; x_0,\wt u^*(\cd)\big)=\xi\,.
\eeq
Proposition \ref{w1225p1} then is used to deduce that
\beq
\|\wt u^*(\cd)\|_{L^2_\dbF(\dbR^+;\dbR^m)}=N(T,x_0,\xi)=N(T(N,x_0,\xi),x_0,\xi)=N\,.
\eeq
Thus, Problem $(TP)_{x_0,\xi}^N$ has a minimal time control.

\ms

{\bf (2)} By using $T=T(N,x_0,\xi)$,
for an arbitrarily fixed minimal time  control $u_*(\cd)$ of Problem $(TP)_{x_0,\xi}^N$, we know that $u_*|_{(0,T)}(\cd)$, the restriction of $u_*(\cd)$ to the interval $(0,T)$, is an admissible control of Problem 
$(NP)_{x_0,\xi}^T$, {\em i.e.},
\bel{w1225e25}
x\big(T, x_0, u_*|_{(0,T)}(\cd)\big)
=\xi \q\mbox{ and }\q 
\|u_*|_{(0,T)}(\cd)\|_{L^2_\dbF(0,T;\dbR^m)}\leq \|u_*(\cd)\|_{L^2_\dbF(\dbR^+;\dbR^m)}\leq N\,.
\ee
Then
\bel{w1225e26}
N\big(T(N,x_0,\xi),x_0,\xi\big)=N(T,x_0,\xi)\leq \|u_*|_{(0,T)}(\cd)\|_{L^2_\dbF(0,T;\dbR^m)}\leq N\,,
\ee
which, together with \rf{w1225e22} in Proposition \ref{w1225p1}, yields
\bel{w1225e27}
N(T,x_0,\xi)= \|u_*|_{(0,T)}(\cd)\|_{L^2_\dbF(0,T;\dbR^m)}= N\,.
\ee
Hence, $u_*|_{(0,T)}(\cd)$ is a minimal norm contro of Problem $(NP)^T_{x_0,\xi}$.

{\bf (3)}
For the minimal norm control $u^*(\cd)$ of Problem $(NP)_{x_0, \xi}^T$, denote by $\wt u^*(\cd)$ the zero extension of $u^*(\cd)$ 
on $[0,\infty)$. Then $\wt u^*(\cd)$ is an admissible control of Problem $(TP)_{x_0,\xi}^{N(T,x_0,\xi)}$, {\em i.e.},
\bel{w1225e28}
x\big(T,x_0,\wt u^*(\cd)\big)=\xi \q\mbox{ and }\q \|\wt u^*(\cd)\|_{L^2_\dbF(\dbR^+;\dbR^m)}=\|u^*(\cd)\|_{L^2_\dbF(0,T;\dbR^m)}=N(T,x_0,\xi)\,.
\ee
By the second conclusion of \rf{w1225e28} and \rf{w1225e22}, we know that 
\beq
\|\wt u^*(\cd)\|_{L^2_\dbF(\dbR^+;\dbR^m)}=N(T,x_0,\xi)
=N(T(N,x_0,\xi),x_0,\xi)=N\,,
\eeq
which yields that $\wt u^*(\cd)$ is a minimal time  control of Problem $(TP)_{x_0,\xi}^N$.
That completes the proof.
\end{proof}

\section{Equivalence of minimal time and norm control problems with zero target}\label{zero-target}

In this section, we consider simplified versions of Problems $(TP)_{x_0,\xi}^N$ and $(NP)_{x_0,\xi}^T$ with target $\xi=0$. 
To shorten notations, we will respectively write  $(TP)_{x_0}^N$, $(NP)_{x_0}^T$, 
$T(N,x_0)$ and $N(T,x_0)$, instead of $(TP)_{x_0,0}^N$, $(NP)_{x_0,0}^T$, 
$T(N,x_0,0)$ and $N(T,x_0,0)$.

When $x_0=0$, Lemma \ref{w1213l1} implies 
that for any $T\in \dbR^+$, the corresponding minimal norm control $u^*(\cd)=0$, and 
consequently $N(T, 0)=0$.
Hence, we focus only on the nontrivial case where $x_0\neq 0$.
The following result concerns the regularity of $N(\cd,x_0)$.

\bl{reg-N}
Suppose that the assumption {\rm \bf (A)} holds, and $x_0\neq 0$. Then
\begin{enumerate}

\item [{\rm (i)}] $N(\cd,x_0):\,\dbR^+\to \dbR^+$ is decreasing;


\item [{\rm (ii)}] $\ds \lim_{T\to 0}N(T,x_0)=+\infty\,.$
\end{enumerate}
\el

\begin{proof}
{\bf (1) Verification of (i).} 
By Lemma \ref{w1213l1}, we know that for any $T>0$,  the minimal norm $N(T,x_0)=\|u^*(\cd)\|_{L^2_\dbF(0,T;\dbR^m)}>0$.

%
%

Suppose that $T_1,T_2\in \dbR^+$, $T_1<T_2$ and $u_1(\cd)$ resp.~$u_2(\cd)$ are minimal norm controls of Problems 
$(NP)_{x_0}^{T_1}$ resp.~$(NP)_{x_0}^{T_2}$.
By setting
\beq
\bar u(t)=
\lt\{\!\!\!
\begin{array}{ll}
u_1(t)\qq & t\in [0,T_1]\,,\\
 0\qq & t\in (T_1,T_2]\,,
\end{array}
\rt.
\eeq
we know that $\bar u(\cd)$ is an admissible control to Problem $(NP)_{x_0}^{T_2}$, and
\beq
N(T_2,x_0)\leq \|\bar u(\cd)\|_{L^2_\dbF(0,T_2;\dbR^m)}=\| u_1(\cd)\|_{L^2_\dbF(0,T_1;\dbR^m)}=N(T_1,x_0)\,.
\eeq
That proves the assertion (i).


\ms

{\bf (2) Verification of (ii).}  By contradiction, suppose that there exists a sequence $\{T_k\}_{k=1}^\infty\subset (0,1)$, $T_k \downarrow 0$, and 
\beq
\lim_{T_k\to0} N(T_k,x_0)=\h N\in (0,+\infty)\,.
\eeq
For any $k$,
suppose that $u_k(\cd)$ is the minimal norm control of Problem $(NP)_{x_0}^{T_k}$. 
By the monotonicity of $N(\cd,x_0)$, we have
\beq
\|u_k(\cd)\|_{L^2_\dbF(0,T_k;\dbR^m)}=N(T_k,x_0)\leq \h N\,.
\eeq

For any $\eta\in L^2_{\mf_{T_1}}(\Omega;\dbR^n)$, by applying It\^o's formula to $\big\lan x\big(t; x_0,u_k(\cd)\big), y(t;\eta, T_1)\big\ran$ and setting 
\bel{w102e1}
\Phi(\cd;\eta,T_1)=B^\top y(\cd;\eta,T_1)+D^\top Y(\cd;\eta,T_1)\,.
\ee
we arrive at
\bel{w727e7}
\lan x_0, y(0;\eta,T_1)\ran +\me\int_0^{T_k}\lan u_k(t), \Phi(t;\eta,T_1)\ran \rd t=0\,.
\ee
Utilizing H\"older's inequality, the fact that $y(\cd;\eta,T_1),Y(\cd;\eta,T_1)\in L^2_\dbF(0,T_1;\dbR^n)$ 
(which implies that $\Phi(\cd;\eta,T_1)$ defined in \rf{w102e1} belongs to $ L^2_\dbF(0,T_1;\dbR^m)$) and 
$T_k\downarrow 0$, we can deduce that
\bel{w727e8}
\bal
&\me\int_0^{T_k}\lan u_k(t), \Phi(t;\eta,T_1)\ran \rd t
\leq \|u_k(\cd)\|_{L^2_\dbF(0,T_k;\dbR^m)} \|\Phi(\cd;\eta,T_1)\|_{L^2_\dbF(0,T_k;\dbR^m)}\\
&\qq\leq \h N   \|\Phi(\cd;\eta,T_1)\|_{L^2_\dbF(0,T_k;\dbR^m)}
\to 0\,.
\eal
\ee
Hence, \rf{w727e7} and \rf{w727e8} imply that
\beq
\lan x_0, y(0;\eta,T_1)\ran=0 \qq \forall\, \eta\in L^2_{\mf_{T_1}}(\Omega;\dbR^n)\,,
\eeq 
which leads to $x_0=0$ contradicting the assumption $x_0\neq 0$,  and then settles the assertion (ii).
\end{proof}

\br{w1118r1}
{\rm (1)}
Under the assumption {\rm\bf (A)},  for any $x_0\neq 0$,  by Lemma \ref{reg-N} {\rm (i)} and the nonnegativity of $N(\cd, x_0)$,  we can define
\beq
\h N(x_0)=\lim_{T\to \infty }N(T,x_0)\,.
\eeq
Furthermore, Theorem \ref{w1230t1} and Lemma \ref{reg-N} imply that
$N(\cd, x_0):\, \dbR^+\to [\h N(x_0),\infty)$ is continuous. 

{\rm (2)}  When $\xi\neq 0$, the decreasing property of $N(\cd; x_0,\xi)$ cannot be established, as demonstrated by the following example.
For the system \rf{sys-1} with 
\beq
A=1,\, B=(0, 1),\, C=0,\,D=(1,0)\,, x_0=-\frac 1 2\,,\xi=-1\,,
\eeq
we can follow the procedure outlined in Remark \ref{w1229r1} to obtain
\beq
\Sigma(t)=\frac {1-\exp\{-2(T-t)\}}{2}\,,\,\, \big(\f(t),\b(t)\big)=\big(\exp\{-(T-t)\}, 0\big)\,, \,\, \Pi(t)= \exp\{-t\}\,,\q \forall\,t\in [0,T]\,,
\eeq
and subsequently,
\beq
\l^*= \frac{2}{\exp\{-2T\}-1}\big(\exp\{-T\}-\frac 1 2\big)\,, \,\, y(t)=\exp\{-t\}\l^*\,,\qq \forall\,t\in[0,T]\,.
\eeq
The optimal norm control of Problem $(NP)^T_{x_0,\xi}$ can then be derived by \rf{w1229e13} as
\beq
u^*(t)=\big(z^*(t),v^*(t)\big)=\big(0, \exp\{-t\}\l^*\big)\,,\qq \forall\,t\in [0,T]\,.
\eeq
Thus, the minimal norm function is given by
\beq
N(T,x_0,\xi)= \sqrt{\frac{2}{1-\exp\{-2T\}}} \Big|\exp\{-T\}-\frac 1 2\Big|\,,
\eeq
from which we observe that $N(\cd,x_0,\xi)$ is decreasing on $(0, \ln 2)$ and then incresing on $(\ln2, \infty)$.

\er

Under the setting $\xi=0$,
the following result improves the result in Proposition \ref{w1225p1}. 

\bp{equivalence-1}
Let the assumption {\rm \bf (A)} hold and $x_0\neq 0$, and suppose that $N>\h N(x_0)$ and $0<T<\infty$. 
Then, the following two conditions are equivalent
\begin{enumerate}
\item [{\rm (i)}]
$
N=N(T,x_0)\,;
$

\item [{\rm (ii)}]
$
T=T(N,x_0)\,.
$

\end{enumerate}
\ep

\begin{proof}
{\bf (i)$\Rightarrow$(ii)}. 
Suppose that $u^*(\cd)\in L^2_\dbF(0,T;\dbR^m)$ is the minimal norm control to Problem $(NP)_{x_0}^{T}$ such 
that
\beq
\|u^*(\cd)\|_{L^2_\dbF(0,T;\dbR^m)}=N(T,x_0)=N\qq \mbox{and}\qq
x\big(T;x_0,u^*(\cd)\big)=0\,,
\eeq
which means that
\bel{w728e4}
T(N,x_0)\leq T<\infty\,,
\ee
and hence, Problem $(TP)_{x_0}^{N}$ admits admissible controls.

Set $T_0=T(N,x_0)$.
By the definition of $T(N,x_0)$ in \rf{minimal-time},
there exists 
a decreasing sequence $\{T_k\}_{k=1}^\infty\subset ( T_0, T_0+1)$ and a sequence
$\{u_k(\cd)\}_{k=1}^\infty\subset \cU_N$ such that $x\big(T_k;x_0,u_k(\cd)\big)=0$, which yields that there exists a subsequence of $\{u_k(\cd)\}_{k=1}^\infty$ denoted in the same way satisfying
\beq
u_k(\cd)\stackrel{\mathrm{w}}{\longrightarrow} u_*(\cd) \qq\mbox{in }L^2_\dbF(\dbR^+;\dbR^m)\,,
\eeq
and
\bel{w1226e5}
\|u_*(\cd)\|_{L^2_\dbF(\dbR^+;\dbR^m)}\leq \liminf_{k\to \infty}\|u_k(\cd)\|_{L^2_\dbF(\dbR^+;\dbR^m)}\leq N\,.
\ee
This shows that $u_*(\cd)\in \cU_N$.

For any $\eta \in L^2_{\mf_{T_0}}(\Omega;\dbR^n)\subset L^2_{\mf_{T_k}}(\Omega;\dbR^n)$, applying It\^o's formula to $\big\lan x\big(t;x_0,u_k(\cd)\big), y(t;\eta, T_k)\big\ran$ yields
\bel{w728e2}
0
=\lan x_0,y(0;\eta,T_k)\ran+\me\int_0^{T_k}\lan u_k(t),\Phi(t;\eta,T_k)\ran\rd t\,,
\ee
where $\big(y(\cd;\eta,T_k), Y(\cd;\eta,T_k)\big)$ solves BSDE \rf{bsde-1} with $T=T_k$, and 
$\Phi(\cd;\eta,T_k)=B^\top y(\cd;\eta,T_k)+D^\top Y(\cd;\eta,T_k)$.
To establish the convergence of \rf{w728e2}, we estimate each term on the right-hand side. Specifically, we decompose the integral term
into three components:
\beq
&\me\int_0^{T_k}\lan u_k(t),\Phi(t;\eta,T_k)\ran\rd t\\
&\q=\me\int_0^{T_0}\lan u_k(t),\Phi(t;\eta,T_k)-\Phi(t;\eta,T_0)\ran\rd t
+\me\int_0^{T_0}\lan u_k(t),\Phi(t;\eta,T_0)\ran\rd t\\
&\qq\q
+\me\int_{T_0}^{T_k}\lan u_k(t),\Phi(t;\eta,T_k)\ran\rd t\\
&\q=:\sum_{i=1}^3 I_i\,.
\eeq
By stability and linearity of BSDEs, we have
\bel{w1112e1}
\bal
\|\Phi(\cd;\eta,T_k)-\Phi(\cd;\eta, T_0)\|^2_{L^2_\dbF(0, T_0;\dbR^m)}
&\leq \cC\|y( T_0;\eta, T_k)-y( T_0;\eta,  T_0)\|^2_{L^2_{\mf_{ T}}(\Omega;\dbR^n)}\\
&= \cC\, \me\big[\|y( T_0;\eta, T_k)-\eta\|^2\big]\,.
\eal
\ee
Based on BSDE \rf{bsde-1}, we proceed to get
\beq
& \me\big[\|y( T_0;\eta, T_k)-\eta\|^2\big]\\
&\q= \me\Big[\Big\|-\int_{T_0}^{T_k} \big[A^\top y(t;\eta,T_k)+C^\top Y(t; \eta,T_k)\big]\rd t +\int_{ T_0}^{T_k} Y(t;\eta, T_k)\rd W(t) \Big\|^2\Big]\\
&\q\leq \cC (T_k- T_0)\me\int_0^{T_k}\big[\|y(t;\eta,T_k)\|^2+\|Y(t;\eta,T_k)\|^2\big]\rd t+\cC\, \me\int_{ T_0}^{T_k} \|Y(t;\eta,T_k)\|^2\rd t\,,
\eeq
where $\cC$ is independent of $k$. This estimate, together with \rf{w1112e1} and the following 
\bel{w1112e2}
\|y(\cd;\eta,T_k)\|_{L^2_\dbF(\O; C([0,T_k];\dbR^n))}+\|Y(\cd;\eta,T_k)\|_{L^2_\dbF(0,T_k;\dbR^n)}\leq \cC\|\eta\|_{L^2_{\mf_{T_0}}(\Omega;\dbR^n)}\,,
\ee
where $\cC$ is independent of $k$,
implies that
\beq
\|\Phi(\cd;\eta,T_k)-\Phi(\cd;\eta,T_0)\|^2_{L^2_\dbF(0, T_0;\dbR^m)}
\leq \cC\,\me\big[\|y( T_0;\eta, T_k)-\eta\|^2\big] \to 0\,.
\eeq
Thus, $I_1\to 0$. Similarly, we conclude by \rf{w1112e2} that $I_3\to 0$. The weak convergence of $\{u_k(\cd)\}$ yields 
$I_2\to \me\int_0^{T_0}\lan u_*(t),\Phi(t;\eta,T_0)\ran\rd t$. 
In a similar vein, we can settle $\lan x_0,y(0;\eta,T_k)\ran\to \lan x_0,y(0;\eta,T_0)\ran$. Combining \rf{w728e2} and the above convergence results, we 
have
\bel{w1112e3}
0=\lan x_0,y(0;\eta,T_0)\ran+\me\int_0^{T_0}\lan u_*(t),\Phi(t;\eta,T_0)\ran\rd t\,.
\ee
Additionally, by applying It\^o's formula, we have
\beq
\me\big\lan x\big( T_0;x_0, u_*(\cd)\big), \eta\big\ran
=\lan x_0, y(0;\eta, T_0)\ran +\me\int_0^{ T_0}\lan u_*(t),\Phi(t; \eta, T_0)\ran \rd t\,,
\eeq
which, together with \rf{w1112e3} and the fact that $\eta\in L^2_{\mf_{ T_0}}(\Omega;\dbR^n)$ is arbitrary, leads to 
\bel{w1112e4}
x\big( T_0;x_0, u_*(\cd)\big)=x\big( T(N,x_0);x_0, u_*(\cd)\big)=0\,.
\ee
Hence, a combination of \rf{w1226e5} and \rf{w1112e4} yields that $u_*(\cd)$ is a minimal time control of 
Problem $(TP)_{x_0,\xi}^N$.
Then, the optimality of $u_*(\cd)$ to Problem $(TP)_{x_0,\xi}^N$ implies
\beq
N\big(T(N,x_0),x_0\big)\leq \|u_{*}(\cd)\|_{L^2_\dbF(0,T(N,x_0));\dbR^m}\leq N=N(T,x_0)\,.
\eeq
Since $N(\cd,x_0)$ is decreasing, we have
\bel{w728e5}
T(N,x_0)\geq T\,.
\ee
Hence, the condition (ii) can be derived by combining with \rf{w728e4} and \rf{w728e5}.

\ms
{\bf (ii)$\Rightarrow$(i)}. 
Directly derived by Proposition \ref{w1225p1}.
%
That completes the proof.
\end{proof}

\bc{w1212c1}
Under the assumption {\rm \bf (A)},
for any $x_0\neq 0$, $N(\cd,x_0)$ is strictly decreasing on $(0, T_{sup})$, where
\beq
T_{sup}=\inf\big\{ T>0\,\big|\, N(T,x_0)=\h N(x_0)\big\}\,.
\eeq
Moreover,  the minimal time function $T(\cd, x_0)$ is strictly decreasing and continuous on $(\h N(x_0), \infty)$.
\ec

\begin{proof}
Let $0<T_1<T_2<T_{sup}$, and we know by Lemma \ref{reg-N} that
\beq
N(T_1,x_0)\geq N(T_2,x_0)\,.
\eeq

We make use of contradiction to derive the assertion. Suppose that $N(T_1,x_0)= N(T_2,x_0)$, 
and set
\beq
N_0\deq N(T_1,x_0)= N(T_2,x_0)> \h N(x_0)\,.
\eeq
By Proposition \ref{equivalence-1}, we have
\beq
T_1=T(N_0,x_0)=T_2\,,
\eeq
which contradicts our assumption $T_1<T_2$.

The strict monotonicity of $T(\cd, x_0)$ can then be established by using the strict monotonicity of $N(\cd, x_0)$ and Proposition 
\ref{equivalence-1}.
That completes the proof.
\end{proof}

\br{w1219r1}
{\rm (1)} Suppose that $x_0\neq 0$ and $(N,T)\in \big\{\big(N(T,x_0), T\big)\,\big|\, T\in(0,T_{sup})\big\}$. Based on Theorem \ref{w1225t2} and Proposition 
\ref{equivalence-1}, Problems $(TP)_{x_0}^N$ and $(NP)_{x_0}^T$ are mutually equivalent. 

{\rm (2)}
When $T\in (0,T_{sup})$, if $u_{*1}(\cd)$ and $u_{*2}(\cd)$ are minimal time controls of Problem $(TP)_{x_0}^{N(T,x_0)}$, then
\bel{w1219e4}
u_{*1}|_{(0,T)}(\cd)=u_{*2}|_{(0,T)}(\cd)\,.
\ee
Indeed, by Theorem \ref{w1225t2} we know that $u_{*1}|_{(0,T)}(\cd)\,,u_{*2}|_{(0,T)}(\cd)$ are 
both minimal norm controls of Problem $(NP)_{x_0}^T$.
Then the uniqueness of the minimal norm control of Problem $(NP)_{x_0}^T$ (see Lemma \ref{w727l1}) yields \rf{w1219e4}.

{\rm (3)} For any $N_0\in (\h N(x_0),\infty)$, by Corrolary \ref{w1212c1}, there exists a unique $T_0\in (0, T_{sup})$ such that $N(T_0,x_0)=N_0$. By Lemma \ref{w727l1}, Problem $(NP)^{T_0}_{x_0}$ admits a unique optimal norm
control $u^*(\cd)$. Then, by Step {\rm (3)} in the proof of Theorem \ref{w1225t2}, we know that $\wt u^*(\cd)$, the zero extension of $u^*(\cd)$ 
over $\dbR^+$, is a minimal time control of Problem $(TP)^{N_0}_{x_0}$.

\er

\section*{Declarations}

\no{\bf Competing interest} \\
The author declares that there is no conflict of interest regarding the publication of this paper.


\begin{thebibliography}{10}

\bibitem{Bi-Sun-Xiong20}
{\sc X.~Bi, J.~Sun, and J.~Xiong}, {\em Optimal control for controllable
  stochastic linear systems}, ESAIM Control Optim. Calc. Var., 26 (2020),
  pp.~Paper No. 98, 23.

\bibitem{Buckdahn-Quincampoix-Tessitore06}
{\sc R.~Buckdahn, M.~Quincampoix, and G.~Tessitore}, {\em A characterization of
  approximately controllable linear stochastic differential equations}, in
  Stochastic partial differential equations and applications---{VII}, vol.~245
  of Lect. Notes Pure Appl. Math., Chapman \& Hall/CRC, Boca Raton, FL, 2006,
  pp.~53--60.

\bibitem{Dou-Lv19}
{\sc F.~Dou and Q.~L\"u}, {\em Partial approximate controllability for linear
  stochastic control systems}, SIAM J. Control Optim., 57 (2019),
  pp.~1209--1229.

\bibitem{Gashi15}
{\sc B.~Gashi}, {\em Stochastic minimum-energy control}, Systems Control Lett.,
  85 (2015), pp.~70--76.

\bibitem{Goreac16}
{\sc D.~Goreac, A.~C. Grosu, and E.-P. Rotenstein}, {\em Approximate and
  approximate null-controllability of a class of piecewise linear {M}arkov
  switch systems}, Systems Control Lett., 96 (2016), pp.~118--123.

\bibitem{Li-Sun-Xiong19}
{\sc X.~Li, J.~Sun, and J.~Xiong}, {\em Linear quadratic optimal control
  problems for mean-field backward stochastic differential equations}, Appl.
  Math. Optim., 80 (2019), pp.~223--250.

\bibitem{Lim-Zhou01}
{\sc A.~E.~B. Lim and X.~Y. Zhou}, {\em Linear-quadratic control of backward
  stochastic differential equations}, SIAM J. Control Optim., 40 (2001),
  pp.~450--474.

\bibitem{Liu-Peng10}
{\sc F.~Liu and S.~Peng}, {\em On controllability for stochastic control
  systems when the coefficient is time-variant}, J. Syst. Sci. Complex., 23
  (2010), pp.~270--278.

\bibitem{Lv19}
{\sc Q.~L\"{u}}, {\em Well-posedness of stochastic {R}iccati equations and
  closed-loop solvability for stochastic linear quadratic optimal control
  problems}, J. Differential Equations, 267 (2019), pp.~180--227.

\bibitem{Lv-Yong-Zhang12}
{\sc Q.~L\"u, J.~Yong, and X.~Zhang}, {\em Representation of {I}t\^o integrals
  by {L}ebesgue/{B}ochner integrals}, J. Eur. Math. Soc. (JEMS), 14 (2012),
  pp.~1795--1823.

\bibitem{Lv-Zhang21}
{\sc Q.~L\"{u} and X.~Zhang}, {\em Mathematical control theory for stochastic
  partial differential equations}, vol.~101 of Probability Theory and
  Stochastic Modelling, Springer, Cham, 2021.

\bibitem{Peng94}
{\sc S.~Peng}, {\em Backward stochastic differential equation and exact
  controllability of stochastic control systems}, Progr. Natur. Sci. (English
  Ed.), 4 (1994), pp.~274--284.

\bibitem{Qin-Wang18}
{\sc S.~Qin and G.~Wang}, {\em Equivalence between minimal time and minimal
  norm control problems for the heat equation}, SIAM J. Control Optim., 56
  (2018), pp.~981--1010.

\bibitem{Wang-Xu-Zhang15}
{\sc G.~Wang, Y.~Xu, and Y.~Zhang}, {\em Attainable subspaces and the bang-bang
  property of time optimal controls for heat equations}, SIAM J. Control
  Optim., 53 (2015), pp.~592--621.

\bibitem{Wang-Zuazua12}
{\sc G.~Wang and E.~Zuazua}, {\em On the equivalence of minimal time and
  minimal norm controls for internally controlled heat equations}, SIAM J.
  Control Optim., 50 (2012), pp.~2938--2958.

\bibitem{Wang-Yang-Yong-Yu17}
{\sc Y.~Wang, D.~Yang, J.~Yong, and Z.~Yu}, {\em Exact controllability of
  linear stochastic differential equations and related problems}, Math. Control
  Relat. Fields, 7 (2017), pp.~305--345.

\bibitem{Wang-Yu20}
{\sc Y.~Wang and Z.~Yu}, {\em On the partial controllability of {SDE}s and the
  exact controllability of {FBSDES}}, ESAIM Control Optim. Calc. Var., 26
  (2020), pp.~Paper No. 68, 27.

\bibitem{Wang-Zhang15}
{\sc Y.~Wang and C.~Zhang}, {\em The norm optimal control problem for
  stochastic linear control systems}, ESAIM Control Optim. Calc. Var., 21
  (2015), pp.~399--413.

\bibitem{Yong-Zhou99}
{\sc J.~Yong and X.~Y. Zhou}, {\em Stochastic controls: Hamiltonian systems and
  HJB equations}, vol.~43 of Applications of Mathematics (New York),
  Springer-Verlag, New York, 1999.

\end{thebibliography}

\end{document}